\documentclass[a4paper,10pt,reqno]{amsart}
\usepackage{a4wide}
\usepackage[english]{babel}
\usepackage{mathrsfs}
\usepackage{amsmath,amssymb,amsthm}
\usepackage{thmtools}
\usepackage{upref}

\usepackage{hyperref}
\usepackage{graphicx}
\usepackage{subcaption}
\usepackage[sort]{natbib}
\let\cite\citep
\usepackage{listings}
\usepackage{xcolor}

\definecolor{codegreen}{rgb}{0,0.6,0}
\definecolor{codegray}{rgb}{0.5,0.5,0.5}
\definecolor{codepurple}{rgb}{0.58,0,0.82}
\definecolor{backcolour}{rgb}{0.95,0.95,0.92}

\lstdefinestyle{mystyle}{
language=Python,
backgroundcolor=\color{backcolour},
commentstyle=\color{codegreen},
keywordstyle=\color{magenta},
numberstyle=\tiny\color{codegray},
stringstyle=\color{codepurple},
basicstyle=\ttfamily\footnotesize,
breakatwhitespace=false,
breaklines=true,
captionpos=b,
keepspaces=true,
numbers=left,
numbersep=5pt,
showspaces=false,
showstringspaces=false,
showtabs=false,
tabsize=2
}
\declaretheorem[numbered=no,name=Barabanov's Theorem]{barthm}
\declaretheorem[numbered=no,name=Main Claim]{mainclaim}
\newtheorem{lemma}{Lemma}
\newtheorem{question}{Question}

\newcommand{\setA}{\mathscr{A}}

\title[Non-Sturmian sequences of matrices]{Non-Sturmian sequences of matrices
providing the maximum growth rate of matrix products}

\author{Victor Kozyakin}

\address{Institute for Information Transmission Problems, Russian Academy of
Sciences, Bolshoj Karetny lane 19, Moscow 127051, Russia}

\email{kozyakin@iitp.ru}



\date{}

\begin{document}
\begin{abstract}
In the theory of linear switching systems with discrete time, as in other
areas of mathematics, the problem of studying the growth rate of the norms of
all possible matrix products $A_{\sigma_{n}}\cdots A_{\sigma_{0}}$ with
factors from a set of matrices $\mathscr{A}$ arises. So far, only for a
relatively small number of classes of matrices $\mathscr{A}$ has it been possible
to accurately describe the sequences of matrices that guarantee the maximum
rate of increase of the corresponding norms. Moreover, in almost all cases
studied theoretically, the index sequences $\{\sigma_{n}\}$ of matrices
maximizing the norms of the corresponding matrix products have been shown
to be periodic or so-called Sturmian, which entails a whole set of ``good''
properties of the sequences $\{A_{\sigma_{n}}\}$, in particular the
existence of a limiting frequency of occurrence of each matrix factor
$A_{i}\in\mathscr{A}$ in them. In the paper it is shown that this is not
always the case: a class of matrices is defined consisting of two $2\times 2$
matrices, similar to rotations in the plane, in which the sequence
$\{A_{\sigma_{n}}\}$ maximizing the growth rate of the norms
$\|A_{\sigma_{n}}\cdots A_{\sigma_{0}}\|$ is not Sturmian. All
considerations are based on numerical modeling and cannot be considered
mathematically rigorous in this part; rather, they should be interpreted as
a set of questions for further comprehensive theoretical analysis.
\end{abstract}

\maketitle

\section{Introduction}\label{S:intro}

Various problems of mathematics reduce to the problem of computing the
maximum growth rate of the norms of matrix products $A_{\sigma_{n}}\cdots
A_{\sigma_{0}}$ with factors from a set of matrices~$\setA$.

One of the basic, though greatly simplified, examples of this type of
situation is found in systems and control theory~\cite{BrayTong:TCS79,
BT:Autom00, BC:TCS03, SWMWK:SIAMREV07, Jungers:09, WuHe:SIAMJCO20} when
considering the asymptotic behavior of solutions of the so-called linear
switching system with discrete time, whose dynamics is described by the
equation
\begin{equation}\label{E:switch} x_{n+1}=A_{\sigma_{n}}x_{n},\qquad
\sigma_{n}\in\{0,1,\ldots,m-1\},~n\ge 0,
\end{equation}
where $A_{\sigma_{i}}\in\setA:=\{A_{0},A_{1},\ldots,A_{m-1}\}$. The solutions
for system~\eqref{E:switch} may be represented as follows:
\begin{equation}\label{E:traj}
x_{n}=A_{\sigma_{n-1}}\cdots A_{\sigma_{1}}A_{\sigma_{0}}x_{0}.
\end{equation}
Therefore, in studying the question of their asymptotic behavior, we
naturally come to the problem of estimating (and preferably computing
exactly) the growth rate of the norms $\|A_{\sigma_{n}}\cdots
A_{\sigma_{1}}A_{\sigma_{0}}\|$ with arbitrary factors
$A_{\sigma_{i}}\in\setA$. Of course, the set of questions related to the
analysis of the asymptotic behavior of $x_{n}$ elements can be extended, but
in this paper we will not deal with such generalizations.

It is worth noting that the problem of computing the maximum possible growth
rate of the norms of matrix products with factors from a set of matrices is
quite general; in particular, numerous problems in other areas of science are
reduced to it, for example, in coding theory~\cite{MOS:IEEETIT01,
BJP:IEEETIT06}, computational mathematics~\cite{DaubLag:LAA92, HStr:LASP92,
Maesumi:AT98, DaubLag:LAA01, JPB:LAA08}, and the theory of parallel and
distributed computation~\cite{ChM:LAA69, BertTsi:89}.

Currently, the range of questions related to the analysis of the growth rate
of the norms of matrix products $A_{\sigma_{n}}\cdots
A_{\sigma_{1}}A_{\sigma_{0}}$ is usually considered in the framework of the
so-called theory of joint/generalized spectral radius, which emerged in the
1960s~\cite{RotaStr:IM60, DaubLag:LAA92, LagWang:LAA95, DaubLag:LAA01} and
now has several hundred publications~\cite{Koz:IITP13}. Moreover, in almost
all cases studied theoretically, the index sequences $\{\sigma_{n}\}$ of
matrices maximizing the norms of the corresponding matrix products turned out
to be periodic or so-called Sturmian. Both the periodicity and the fact that
the index sequences $\{\sigma_{n}\}$ are Sturmian entail a whole set of
``good'' properties of the sequences $\{A_{\sigma_{n}}\}$, in particular the
existence of a limiting frequency of occurrence of each matrix factor
$A_{i}\in\setA$~\cite{BM:JAMS02, BTV:SIAMJMA03, Koz:CDC05:e,
Koz:INFOPROC06:e}.

This may give the false impression that periodic or Sturmian sequences occur
every time one tries to maximize the norms of matrix products (at least in
the case of a pair of ${2\times2}$ matrices). This impression is reinforced
by the fact that we are not aware of any theoretical studies in this area,
apart from those leading to the appearance of periodic or Sturmian sequences,
which can be explained by the extreme theoretical and technical complexity of
the corresponding analysis. The aim of this paper is to refute this
impression. To this end, we determine a class of ${2\times2}$ matrices
consisting of two matrices similar to rotations of the plane in which the
sequence $\{A_{\sigma_{n}}\}$ maximizing the growth rate of the norms
$\|A_{\sigma_{n}}\cdots A_{\sigma_{0}}\|$ is not Sturmian.

Let us describe the structure of the work. The introduction provides a
rationale for the topics covered here. Section~\ref{S:background} briefly
recalls the main facts and constructions from the theory of joint/generalized
spectral radius, among which \nameref{T:Bar} plays a crucial role.
Section~\ref{S:barnorm} recalls the concept of extremal trajectories, i.e.,
trajectories with the maximum rate of increase in a certain Barabanov norm.
Here we describe a general approach which, in the case of a pair of
${2\times2}$ matrices, reduces the problem of constructing extremal
trajectories to the problem of studying iterations of a certain mapping of an
interval into itself (or, equivalently, of a circle into itself) with a
fairly simple structure. Section~\ref{S:sturmian-case} recalls the well-known
theoretical results on the construction and growth of extremal trajectories,
which refer to the case of a pair of nonnegative ${2\times2}$ matrices of a
special form. This case underlies most modern studies of ``nontrivial''
situations in joint/generalized spectral radius theory. The above theoretical
results are illustrated by examples of computer simulations.
Section~\ref{S:rotmaps} considers the case of a pair of matrices, each of
which is similar to a rotation matrix. Using the results of numerical
simulations, it is shown that for such matrices a previously unobserved
phenomenon occurs in which the index sequences of the extremal trajectories
turn out to be non-Sturmian (\nameref{P:main}). For a description of the
methods and means for numerical modeling of the behavior of matrix products
used in this work, see Section~\ref{S:numeric}.

\section{Theoretical Background}\label{S:background}

Recall the basic concepts and results related to the theory of
joint/generalized spectral radius, following~\cite{Koz:CDC05:e,
Koz:INFOPROC05:e, Koz:INFOPROC06:e}.

Let $\setA=\{A_{0},\ldots,A_{m-1}\}$ be a set of $m$ real $d\times d$
matrices, and $\|\cdot\|$ be some norm in ${\mathbb{R}}^{d}$. For every
$n\ge1$, with every finite sequence $\boldsymbol{\sigma}=
\{\sigma_{0},\sigma_{1},\ldots,\sigma_{n-1}\}\in{\{0,\ldots,m-1\}}^{n}$ we
link the matrix
\[
A_{\boldsymbol{\sigma}}=A_{\sigma_{n-1}}\cdots A_{\sigma_{1}}A_{\sigma_{0}},
\]
and define two numerical values:
\[
\rho_{n}({\setA})=\max_{\boldsymbol{\sigma}\in{\{0,\ldots,m-1\}}^{n}}
\|A_{\boldsymbol{\sigma}}\|^{1/n},\qquad
\bar{\rho}_{n}({\setA})=\max_{\boldsymbol{\sigma}\in{\{0,\ldots,m-1\}}^{n}}
{\rho(A_{\boldsymbol{\sigma}})}^{1/n},
\]
where $\rho(\cdot)$ denotes the spectral radius of a matrix. In these
designations, the limit
\[
\rho({\setA})= \limsup_{n\to\infty}\rho_{n}({\setA}),
\]
which does not depend on the choice of the norm $\|\cdot\|$ (and in fact
coincides with the limit $\rho({\setA})=\lim_{n\to\infty}\rho_{n}({\setA})$)
is called the \emph{joint spectral radius} of the set of matrices
$\setA$~\cite{RotaStr:IM60}. Similarly, we can consider the limit
\[
\bar{\rho}({\setA})= \limsup_{n\to\infty}\bar{\rho}_{n}({\setA}),
\]
called the \emph{generalized spectral radius} of the matrix set
$\setA$~\cite{DaubLag:LAA92}. The values $\rho({\setA})$ and
$\bar{\rho}({\setA})$ for bounded families of matrices $\setA$ actually
coincide with each other~\cite{BerWang:LAA92}, and, moreover, for any $n$,
the following inequalities hold:
\begin{equation}\label{Eq-sprad}
\bar{\rho}_{n}({\setA})\le \bar{\rho}({\setA})=\rho({\setA})\le
\rho_{n}({\setA}).
\end{equation}

It follows from the definition of the joint spectral radius that for each
$\varepsilon>0$ the rate of growth of the norms
$\|A_{\boldsymbol{\sigma}}\|=\|A_{\sigma_{n-1}}\cdots
A_{\sigma_{1}}A_{\sigma_{0}}\|$ for large $n$ does not exceed
${(\rho({\setA})+\varepsilon)}^{n}$, that is,
\begin{equation}\label{E:upbound}
\|A_{\boldsymbol{\sigma}}\|=\|A_{\sigma_{n-1}}\cdots
A_{\sigma_{1}}A_{\sigma_{0}}\|\le {(\rho({\setA})+\varepsilon)}^{n}
\end{equation}
for each finite sequence of indices
\[
\boldsymbol{\sigma}=
\{\sigma_{0},\sigma_{1},\ldots,\sigma_{n-1}\}\in{\{0,\ldots,m-1\}}^{n}.
\]
Moreover, there are arbitrarily large $n$ and sequences of indices
$\boldsymbol{\sigma}$ for which\footnote{Here and in the following $\rho(A)$,
where $A$ is a matrix, denotes the spectral radius of this matrix, i.e.\ the
maximum of the absolute values of the eigenvalues of the matrix $A$.}
\[
\rho(A_{\boldsymbol{\sigma}})=\rho(A_{\sigma_{n-1}}\cdots
A_{\sigma_{1}}A_{\sigma_{0}}) \ge {(\bar{\rho}({\setA})-\varepsilon)}^{n},
\]
and therefore, for such $n$, by virtue of~\eqref{Eq-sprad}, the inequalities
\begin{align}\notag
\|A_{\boldsymbol{\sigma}}\|&=\|A_{\sigma_{n-1}}\cdots
A_{\sigma_{1}}A_{\sigma_{0}}\|\ge\\
\label{E:lowbound} &\ge \rho(A_{\sigma_{n-1}}\cdots
A_{\sigma_{1}}A_{\sigma_{0}}) \ge
{(\bar{\rho}({\setA})-\varepsilon)}^{n}\equiv
{(\rho({\setA})-\varepsilon)}^{n}
\end{align}
will hold.

Inequalities~\eqref{E:upbound} and~\eqref{E:lowbound} raise at least two
questions: \emph{first, is it possible to set $\varepsilon$ in them equal to
zero, and second, if this is possible, how can we describe the sets of
indices $\boldsymbol{\sigma}= \{\sigma_{0},\sigma_{1},\ldots,\sigma_{n-1}\}$,
for which inequality~\eqref{E:upbound} becomes an equality with
$\varepsilon=0$?}

The answer to the first question is negative; it follows from the following
simple remark. Let the set $\setA$ consist of one square matrix $A$. Then by
the well-known Gelfand formula~\cite[Corollary~5.6.14]{HJ:e} both values
$\rho({\setA})$ and $\bar{\rho}({\setA})$ coincide with the spectral radius
$\rho(A)$ of the matrix $A$. In this case, by reducing the matrix $A$ to
normal Jordan form, one can easily establish that the growth rate of the
norms $\|A^{n}\|$ is of order ${\rho(A)}^{n}$ if and only if the eigenvalues
of the matrix $A$ that have the largest absolute value are semisimple. At the
same time, for the case when at least one such eigenvalue of the matrix $A$
is not semisimple, the growth rate of the norms $\|A^{n}\|$ is of order
$p(n){\rho(A)}^{n}$, where $p(t)$ is a polynomial.

The answer to the first question has another nuance: Even in cases where the
growth rate of the norms of the matrix products $\|A_{\sigma_{n-1}}\cdots
A_{\sigma_{1}}A_{\sigma_{0}}\|$ for a particular choice of the sequence of
indices
$\{\sigma_{0},\sigma_{1},\ldots,\sigma_{n-1}\}\in{\{0,\ldots,m-1\}}^{n}$ may
coincide with ${\rho(\setA)}^{n}$, this may not happen for all norms of
$\|\cdot\|$ but only for a particular choice of the corresponding norm.
Moreover, even in the case of sets of matrices $\setA$ consisting of a single
matrix, the construction of the corresponding norm turns out to be a
nontrivial problem!

However, even in the simplest nontrivial case, when the set of matrices
$\setA$ consists of a pair of matrices of dimension ${2\times2}$, the
question posed turns out to be much more complicated than the analysis of the
growth rate of the norms of degrees $\|A^{n}\|$ of one matrix $A$. More
precisely, unlike the corresponding analysis for one matrix, in the general
case the computation of the maximal growth rate of the norms
$\|A_{\sigma_{n}}\cdots A_{\sigma_{1}}A_{\sigma_{0}}\|$ turns out to be
algebraically impossible~\cite{Koz:AiT90:6:e, Koz:AiT03:9:e, Koz:ArXiv13},
and the approximate computation of the corresponding rate turns out to be
NP-hard~\cite{TB:MCSS97-1, BT:Autom00}.

Nevertheless, in a fairly general situation, a theoretically satisfactory
answer to the first question can still be obtained. Recall that a set of
matrices $\setA$ is called \emph{irreducible} if matrices from $\setA$ have
no common invariant subspaces other than $\{0\}$ and ${\mathbb{R}}^{d}$. The
irreducibility of a set of matrices plays the same role in studying the
growth rate of the norms of matrix products with multiple matrix factors as
the semisimplicity of the eigenvalues with the largest absolute value of a
single matrix plays in studying the growth rate of the norms of its powers.
The approach proposed by N.~Barabanov in~\cite{Bar:AIT88-2:e, Bar:AIT88-3:e,
Bar:AIT88-5:e} proved to be the most fruitful here, and it has been further
developed in several papers, from which we single out~\cite{Wirth:LAA02}.

\begin{barthm}
\makeatletter\def\@currentlabelname{Barabanov's Theorem}\makeatother%
\label{T:Bar}%
If the set of matrices $\setA=\{A_{0},\ldots,A_{m-1}\}$ is irreducible, then
the number $\rho$ is the joint (generalized) spectral radius of the set
$\setA$ if and only if there exists a norm $\|\cdot\|$ in ${\mathbb{R}}^{d}$
such that
\begin{equation}\label{Eq-mane-bar}
\rho\|x\|= \max\left\{\|A_{0}x\|,\|A_{1}x\|,\ldots,\|A_{m-1}x\|\right\},
\qquad\forall~x\in\mathbb{R}^{d}.
\end{equation}
\end{barthm}
A norm satisfying \eqref{Eq-mane-bar} is usually called the \emph{Barabanov
norm} corresponding to the set of matrices $\setA$. This theorem does not
provide a constructive description of Barabanov norms. Nevertheless, it turns
out to be very effective in analyzing the growth of matrix products. In
particular, it follows from \nameref{T:Bar} that for irreducible sets of
matrices $\setA$, for every finite sequence of indices $\boldsymbol{\sigma}=
\{\sigma_{0},\sigma_{1},\ldots,\sigma_{n-1}\}\in{\{0,\ldots,m-1\}}^{n}$, in
the Barabanov norm the inequality
\[
\|A_{\sigma_{n-1}}\cdots A_{\sigma_{1}}A_{\sigma_{0}}\|\le
{\rho({\setA})}^{n}
\]
holds, and that there is an infinite sequence of indices
$\boldsymbol{\sigma}= \{\sigma_{0},\sigma_{1},\ldots\}$ such that for every
$n$ the equality
\[
\|A_{\sigma_{n-1}}\cdots A_{\sigma_{1}}A_{\sigma_{0}}\|= {\rho({\setA})}^{n}
\]
holds. So for irreducible sets of matrices there are much stronger statements
than inequalities~\eqref{E:upbound} and~\eqref{E:lowbound}.

The second question, \emph{how to describe sets of indices
$\boldsymbol{\sigma}= \{\sigma_{0},\sigma_{1},\ldots,\sigma_{n-1}\}$, for
which inequality~\eqref{E:upbound} turns into equality when $\varepsilon=0$},
is much more complicated than the first. The answer is currently known either
in trivial (and therefore theoretically uninteresting) situations or in some
special and rather restrictive cases. The answer to this question (as of the
current understanding of the problem) is closely related to the concept of
extremal trajectories of sets of matrices, which we briefly recall in the
next section.

\section{Extremal Trajectories}\label{S:barnorm}

In stating the main facts of this section we shall, as in the preceding
section, adhere to the works~\cite{Koz:CDC05:e, Koz:INFOPROC05:e,
Koz:INFOPROC06:e}. Additional discussion of the problems and statements
involved can also be found in~\cite{Jungers:09}.

The question of the growth rate of matrix products with factors from a
certain set of matrices $\setA$ is closely related to a similar question of
the growth rate of solutions of the difference equation~\eqref{E:switch} for
all possible choices of index sequences $\{\sigma_{n}\}$ and initial values
$x_{0}$. The solutions $\{x_{n}\}$ of equation~\eqref{E:switch} will also be
called the \emph{trajectories} defined by the set of matrices $\setA$, or
simply the trajectories of the set of matrices $\setA$. Since, according
to~\eqref{E:traj}, each element $x_{n}$ of the trajectory $\{x_{n}\}$ can be
represented as
\[
x_{n}= A_{\sigma_{n-1}}\cdots A_{\sigma_{1}}A_{\sigma_{0}}x_{0},
\]
then due to~\eqref{E:upbound}
\[
\|x_{n}\|\le \|A_{\sigma_{n-1}}\cdots
A_{\sigma_{1}}A_{\sigma_{0}}\|\,\|x_{0}\|\le
e^{(\rho({\setA})+\varepsilon)n}\|x_{0}\|
\]
for all sufficiently large $n$.

A trajectory $\{x_{n}\}$ of the set of matrices $\setA$ was called
\emph{characteristic} in~\cite{Koz:CDC05:e, Koz:INFOPROC05:e,
Koz:INFOPROC06:e} if it satisfies the inequalities
\[
c_{1}\rho^{n}(\setA)\le\|x_{n}\|\le c_{2}\rho^{n}(\setA), \qquad
n=0,1,2,\ldots\,,
\]
for some choice of constants $c_{1}, c_{2}\in(0,\infty)$. In other words, by
characteristic trajectories are meant those trajectories which, in the
natural sense, are uniformly comparable to the sequence $\{\rho^{n}(\setA)\}$
on the whole infinite interval $n=0,1,2,\ldots$ of the variation of their
indices. Note that the definition of the characteristic trajectory does not
depend on the choice of the norm $\|\cdot\|$ in the space $\mathbb{R}^{d}$.

An important special case of characteristic trajectories is the so-called
extremal trajectories. A trajectory $\{x_{n}\}$ of a collection of matrices
$\setA$ will be called \emph{extremal} if in some Barabanov norm $\|\cdot\|$
it satisfies the identity
\begin{equation}\label{Eq-defextrimseq}
\rho^{-n}(\setA)\|x_{n}\| \equiv \textrm{const}.
\end{equation}

As mentioned in Section~\ref{S:background}, there are Barabanov norms for any
irreducible set of matrices. Then, for any irreducible set of matrices, there
are also extremal trajectories and hence characteristic trajectories. For the
proof, it suffices to construct the trajectory $\{x_{n}\}$ of the set of
matrices $\setA=\{A_{0},\ldots,A_{m-1}\}$, satisfying the initial condition
$x_{0}=x$, recursively. Let the element $x_{n}$ have already been formed.
Then, according to the definition of the Barabanov norm, the equality
\[
\rho(\setA)\|x_{n}\|=
\max\left\{\|A_{0}x_{n}\|,\|A_{1}x_{n}\|,\ldots,\|A_{m-1}x_{n}\|\right\}
\]
holds. Therefore, there exists an index $\sigma_{n}$ such that
\[
\rho(\setA)\|x_{n}\|= \|A_{\sigma_{n}}x_{n}\|,
\]
and for condition~\eqref{Eq-defextrimseq} to be satisfied, it suffices to
define the element $x_{n+1}$ by the equality $x_{n+1}=A_{\sigma_{n}}x_{n}$.

Unlike the definition of a characteristic trajectory, the definition of an
extremal trajectory depends on the choice of an extremal norm: A trajectory
that is extremal in one norm may not be extremal in another norm.
Nevertheless, for an irreducible set of matrices, there are always
\emph{universal extremal trajectories} in a certain sense, i.e., trajectories
which are extremal with respect to any extremal norm: As shown
in~\cite[Th.~3]{Koz:INFOPROC06:e}, every limit point of every normed
characteristic trajectory $\{x_{n}/\|x_{n}\|\}$ serves as the initial value
of a trajectory which is automatically extremal in every Barabanov norm.

The description of extremal trajectories includes, besides the sequence
$\boldsymbol{x}=\{x_{n}\}$, the index sequence $\{\sigma_{n}\}$, which is
used to obtain the trajectory $\{x_{n}\}$ according to~\eqref{E:switch}. In
the following, we describe a construction that allows to define extremal
trajectories as all possible trajectories of a multivalued nonlinear
dynamical system, thus dispensing with the explicit description of the index
sequence $\{\sigma_{n}\}$.

Let $\rho=\rho(\setA)$, and let $\|\cdot\|$ be a Barabanov norm corresponding
to the set of matrices $\setA=\{A_{0},\ldots,A_{m-1}\}$. For each
$x\in\mathbb{R}^{d}$, we define the mapping $g(x)$ by
\begin{equation}\label{Eq-ExtTraject}
g(x):=\{w:~\exists i\in\{0,\ldots,m-1\}:~ w=A_{i}x,\textrm{~where~}
\|A_{i}x\|=\rho\|x\|\}.
\end{equation}
By the definition of the Barabanov norm, for any ${x\in\mathbb{R}^{d}}$, the
set $g(x)$ is nonempty and consists of at most $m$ elements. Note that every
mapping $g(x)$ has a closed graph and the identity
\begin{equation}\label{Eq-genextseq1}
\|g(x)\|\equiv \rho\|x\|
\end{equation} holds for it.

Clearly, the sequence $\boldsymbol{x}=\{x_{n}\}$ is an extremal trajectory of
the set of matrices $\setA$ in the Barabanov norm $\|\cdot\|$ if and only if
it satisfies the inclusions
\[
x_{n+1}\in g(x_{n}),\qquad \forall~n.
\]
In other words, every trajectory of the multivalued mapping $g(\cdot)$ turns
out to be an extremal trajectory of the set of matrices $\setA$ in the
Barabanov norm $\|\cdot\|$. This gives rise to call the map $g(\cdot)$ a
\emph{generator of extremal trajectories}. Just like the Barabanov norm, the
map $g(\cdot)$ cannot be stated explicitly in the general case. Nevertheless,
a rather detailed description of the properties of generators of extremal
trajectories can be obtained for two-element (i.e., $m=2$) sets of
nonnegative ${2\times2}$ matrices. Let us describe the corresponding
constructions in detail.

We fix in space $\mathbb{R}^{2}$ a Barabanov norm $\|\cdot\|$ corresponding
to the set of matrices $\setA=\{A_{0},A_{1}\}$. Let us define the sets
\begin{equation}\label{Eq-defX}
 X_{0}=\{x:~\|A_{0}x\|=\rho\|x\|\},\qquad X_{1}=\{x:~\|A_{1}x\|=\rho\|x\|\}.
\end{equation}
Each of these sets is closed, conic (i.e., it contains, together with the
vector $x\neq0$, each vector of the form $tx$, where $t\ge0$), and by the
definition of the Barabanov norm $X_{0}\cup X_{1}={\mathbb{R}}^{2}$. In this
case, the generator of extremal trajectories $g(\cdot)$
(see~\eqref{Eq-ExtTraject}) in the norm $\|\cdot\|$ takes the form
\begin{equation}\label{Eq-maxmap}
g(x)=\left\{\begin{array}{cl}
  A_{0}x & \textrm{for~}x\in X_{0}\backslash X_{1}, \\
  A_{1}x & \textrm{for~}x\in X_{1}\backslash X_{0},\\
  \{A_{0}x,A_{1}x\}& \textrm{for~}x\in X_{0}\cap X_{1}.
\end{array}  \right.
\end{equation}

Let us examine more closely the structure of the mapping $g(\cdot)$. Let
$(r,\varphi)$ be the polar coordinates of the vector $x\in\mathbb{R}^{2}$. We
denote by $\Omega_{0}$ and $\Omega_{1}$ the angular projections of the conic
sets $X_{0}$ and $X_{1}$, respectively. As mentioned above
(see~\eqref{Eq-genextseq1}), the mapping $g(\cdot)$ satisfies the identity
$\|g(x)\|\equiv \|x\|$. Thus, in the polar coordinate system $(r,\varphi)$,
the mapping $g$ has the form of a mapping with separable variables
\begin{equation}\label{Eq-maxmanlx}
g:(r,\varphi)\mapsto (\rho r,\varPhi(\varphi)),
\end{equation}
where $\rho=\rho(\setA)$ and
\begin{equation}\label{Eq-lxdef}
\varPhi(\varphi)=\left\{\begin{array}{cl}
  \varPhi_{0}(\varphi) & \textrm{for~}\varphi\in \Omega_{0},\\
  \varPhi_{1}(\varphi) & \textrm{for~}\varphi\in \Omega_{1},\\
  \{\varPhi_{0}(\theta),\varPhi_{1}(\theta)\} & \textrm{for~}\varphi\in
  \Omega_{0}\cap \Omega_{1}.
\end{array}  \right.
\end{equation}
Here the functions $\varPhi_{0}(\varphi)$ and $\varPhi_{1}(\varphi)$ are
explicitly defined as the angular coordinates of the mappings $A_{0}x$ and
$A_{1}x$, when $x= (r,\varphi)$ in polar coordinates.

The angular coordinate $\varphi$ characterizes the direction of the vector
$x=(r,\varphi)$. Accordingly, it is natural to interpret $\varPhi(\varphi)$,
$\varphi\in[0,2\pi)$, as a \emph{direction function} or \emph{angular
function} of the generator of extremal trajectories $g(\cdot)$. Note also
that the function $\varPhi(\varphi)$, while $2\pi$-periodic, is in general
not continuous. And since it is obtained as a result of taking the angular
coordinates of linear mappings, its value modulo $\pi$, the function
\begin{equation}\label{E:angfun}
\tilde{\varPhi}(\varphi)=\varPhi(\varphi)\bmod\pi,\qquad\varphi\in[0,\pi),
\end{equation}
is a $\pi$-periodic function.

\sloppy From the definition~\eqref{Eq-ExtTraject} of the function $g(\cdot)$
and its representation in the form~\eqref{Eq-maxmanlx} we obtain the
following description of extremal
trajectories~\cite[Lemma~6]{Koz:INFOPROC06:e}.

\fussy\begin{lemma}\label{L:Phimap} The nonzero trajectory $\{x_{n}\}$ is
extremal for the set of ${{2\times2}}$ matrices $\setA=\{A_{0},A_{1}\}$ in
the Barabanov norm $\|\cdot\|$ if and only if its elements in the polar
coordinate system $(r,\varphi)$ are representable as $x_{n}=(\rho^{n}
r_{0},\varphi_{n})$, where $\rho$ is the joint/generalized spectral radius of
the set of matrices $\setA$, and $\{\varphi_{n}\}$ is the trajectory of the
multivalued mapping $\varPhi(\cdot)$, i.e.
\[
\varphi_{n+1}\in \varPhi(\varphi_{n}),\qquad n=0,1,\ldots~.
\]

Moreover, the trajectory $\{x_{n}\}$ satisfies the equations
\[
x_{n+1}=A_{\sigma_{n}}x_{n},\qquad n=0,1,\ldots~,
\]
with some index sequence $\{\sigma_{n}\}$ if and only if the trajectory
$\{\varphi_{n}\}$ satisfies the equations
\[
\varphi_{n+1}=\varPhi_{\sigma_{n}}(\varphi_{n}),\qquad n=0,1,\ldots~,
\]
or, which is equivalent, when the trajectory $\{\tilde{\varphi}_{n}\}$ with
elements
\[
\tilde{\varphi}_{n}=\varphi_{n}\bmod\pi,\qquad n=0,1,\ldots~,
\]
satisfies the equations
\[
\tilde{\varphi}_{n+1}=\tilde{\varPhi}_{\sigma_{n}}(\tilde{\varphi}),\qquad
n=0,1,\ldots~.
\]
\end{lemma}

\section{A Pair of Nonnegative Matrices}\label{S:sturmian-case}
Despite the fact that the Barabanov norm $\|\cdot\|$ is, as a rule, not known
explicitly, the angular function $\varPhi(\varphi)$ of the generator of
extremal trajectories $g(\cdot)$ turns out in some cases to be defined
``unambiguous enough.'' In this context, we recall some constructions and
results from~\cite{Koz:INFOPROC05:e, Koz:INFOPROC06:e} developed to construct
one of the counterexamples to the so-called Finiteness
Conjecture~\cite{LagWang:LAA95, BM:JAMS02, BTV:SIAMJMA03}.

Consider the set of matrices $\setA=\{A_{0},A_{1}\}$, where
\begin{equation}\label{Eq-M1}
A_{0}=\alpha\left[\begin{array}{cc}
  a & b \\
  0 & 1
\end{array}\right],\quad
A_{1}=\beta\left[\begin{array}{cc}
  1 & 0 \\
  c & d
\end{array}\right].
\end{equation}
For this set of matrices in~\cite{Koz:INFOPROC05:e, Koz:INFOPROC06:e} it was
possible to perform a detailed analysis of the structure of the extremal
trajectories under the additional assumptions that $\alpha,\beta > 0$ and
\[
bc\ge 1\ge a,d > 0.
\]

\begin{figure}[htbp!]
\centering \mbox{}\hfill\subcaptionbox{A Barabanov norm\label{F:0a}}
{\includegraphics*[height=0.31\textwidth]{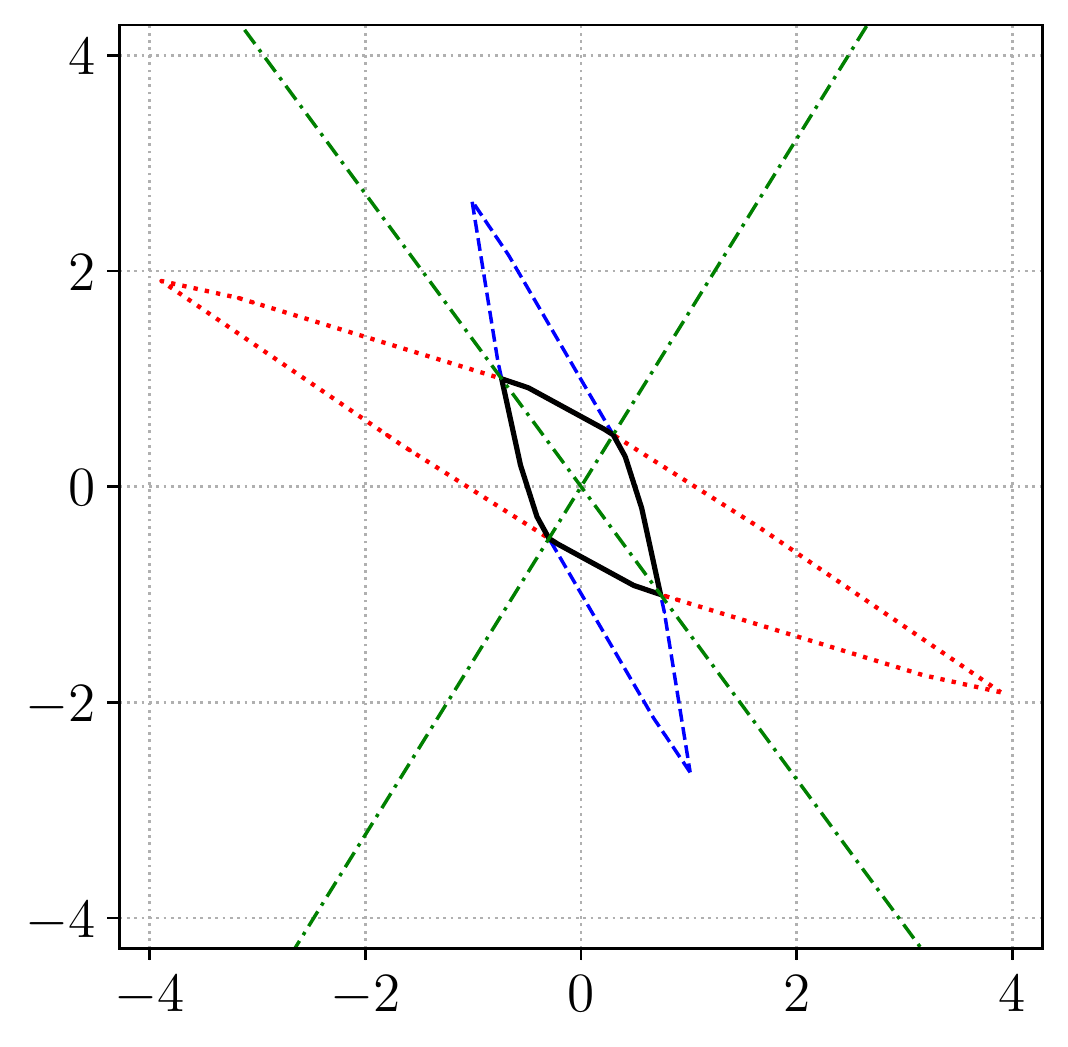}}
\hfill\subcaptionbox{A maximum growth rate trajectory\label{F:0b}}
{\includegraphics*[height=0.31\textwidth]{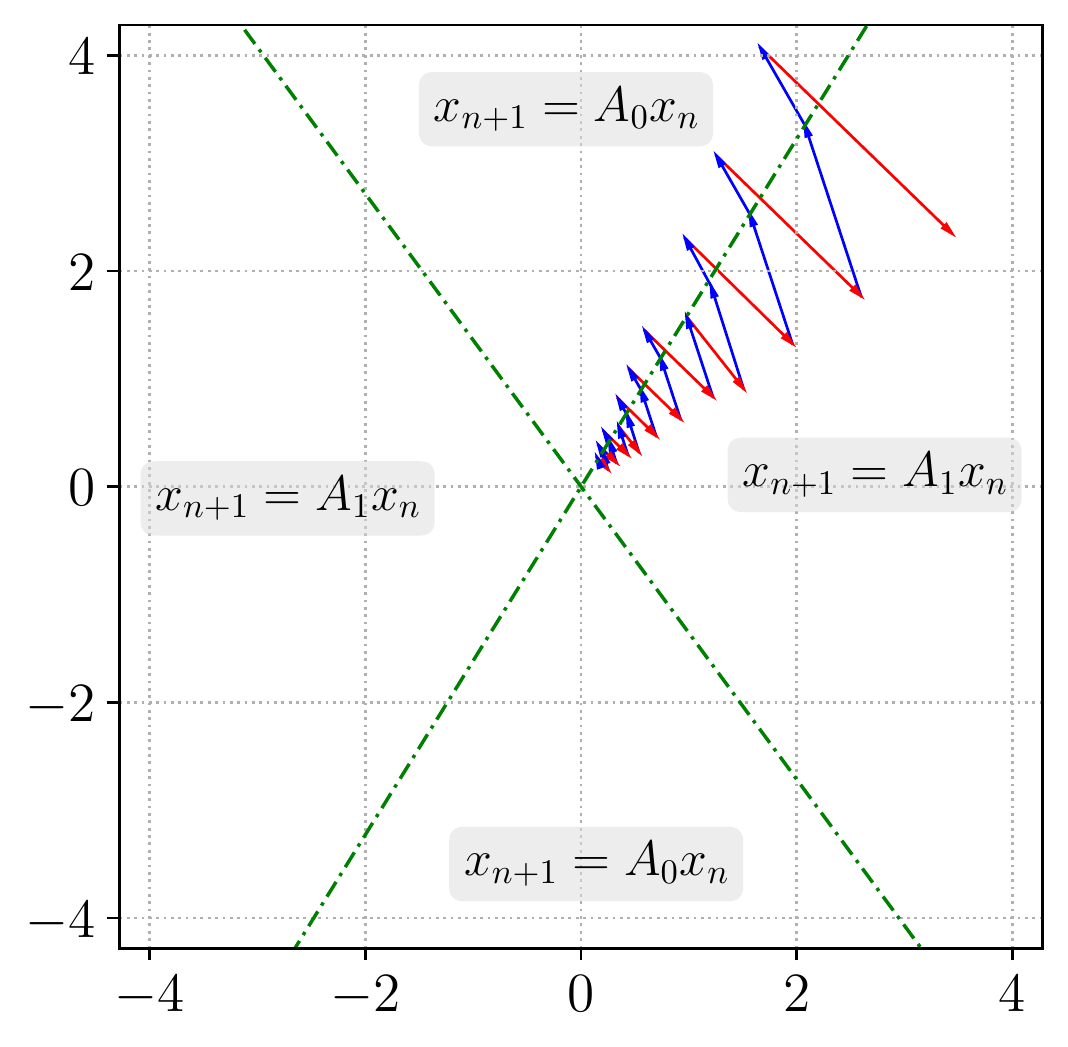}}
\hfill\subcaptionbox{The function $\tilde{\varPhi}(\varphi)$ preserves
orientation and has two discontinuity points\label{F:0c}}
{\includegraphics*[height=0.312\textwidth]{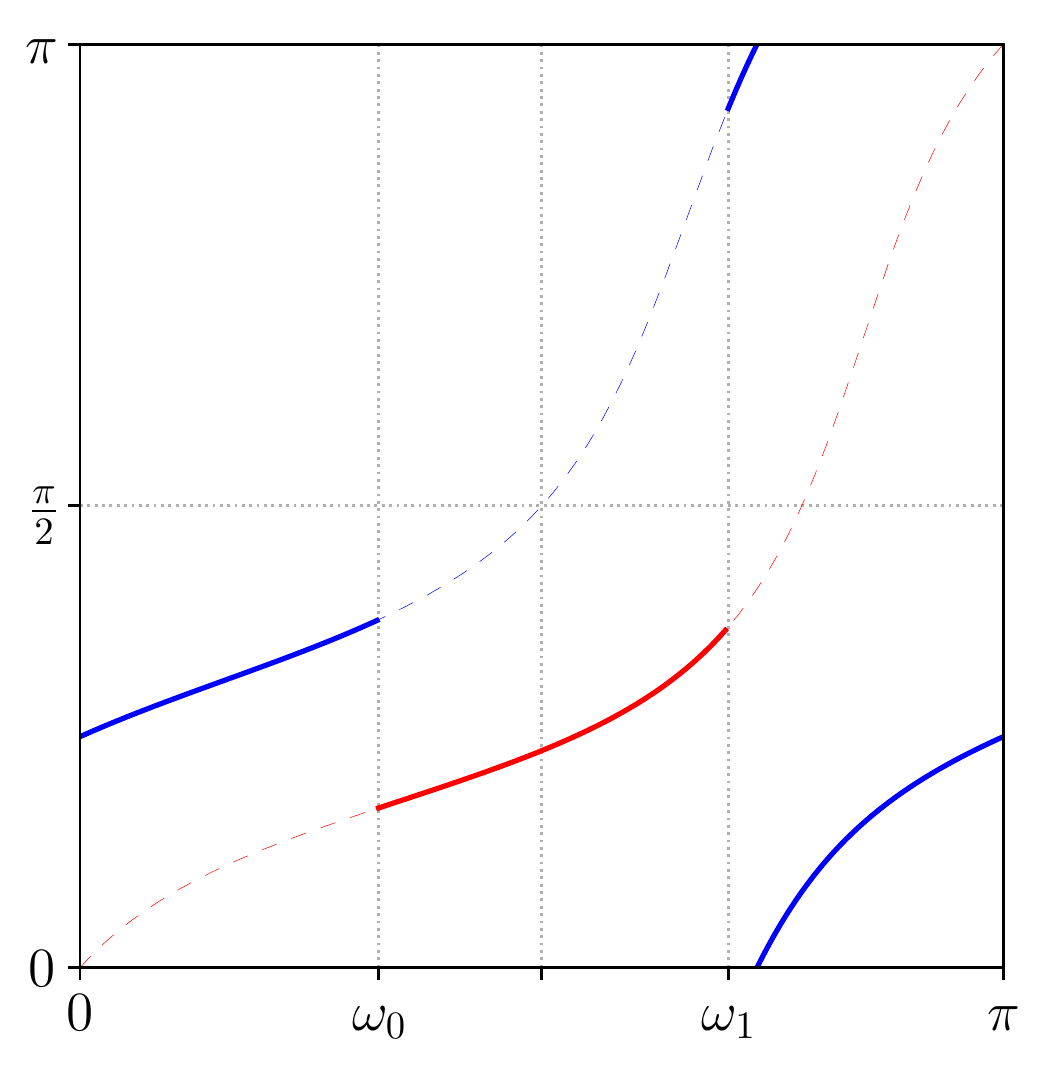}}
\hfill\mbox{} \caption{A Barabanov norm and the angular function for the set
of matrices~\eqref{Eq-M1}, where $\alpha=0.576$, $\beta=0.8$, $a=d=0.9$,
$b=1.1$, $c=1$}\label{F:0}
\end{figure}

\begin{figure}[htbp!]
\centering
\includegraphics*[height=0.31\textwidth]{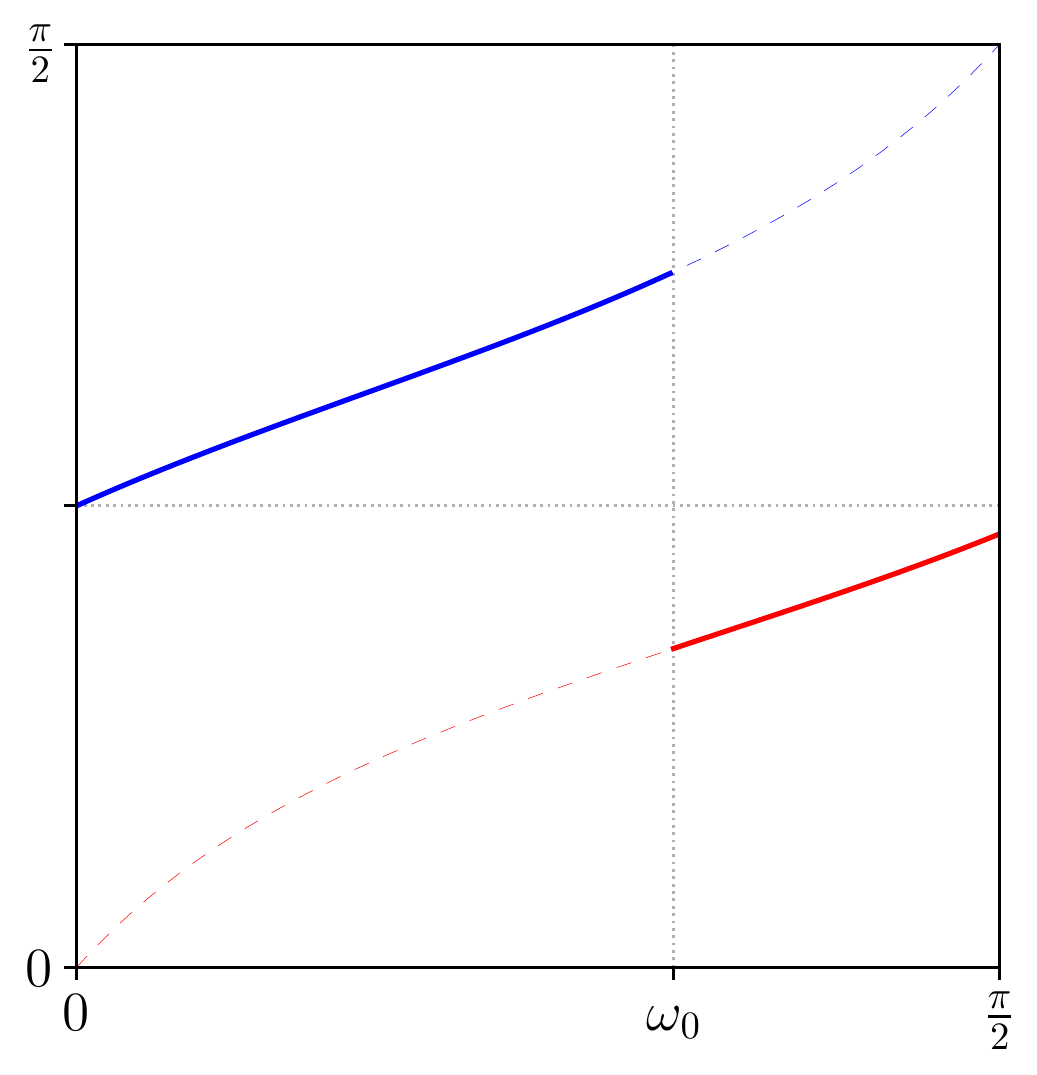}
\caption{The angular function $\tilde{\varPhi}(\varphi)$,
$\varphi\in\left[0,\frac{\pi}{2}\right)$, for the set of
matrices~\eqref{Eq-M1}}\label{F:01}
\end{figure}

In this paper, the approximate construction of the Barabanov norms of the
sets of the matrix sets~\eqref{Eq-M1} (and the visualization of their unit
spheres) were carried out using the algorithms and programs described in
Section~\ref{S:numeric}. An example of the unit sphere of the Barabanov norm
for the set of matrices~\eqref{Eq-M1}, one of the extremal trajectories, and
the corresponding angular function $\tilde{\varPhi}(\varphi)$ for the case
\begin{equation}\label{mycase}
	\alpha=0.576,\quad \beta=0.8,\quad a=d=0.9,\quad b=1.1,\quad c=1,
\end{equation}
is shown in Fig.~\ref{F:0a}. Here the black solid line represents the unit
sphere of the Barabanov norm, i.e., the set of points $x\in\mathbb{R}^{2}$
for which $\|x\|=1$ holds. Dotted and dashed lines denote the sets of points
$x\in\mathbb{R}^{2}$ for which $\|A_{0}x\|=\rho$ and $\|A_{1}x \|=\rho$,
respectively, where $\rho=\rho(\setA)$ is the joint/generalized spectral
radius of the matrix set $\setA$. Dash-dotted lines denote straight lines
consisting of points $x\in X_{0}\cap X_{1}$ (see~\eqref{Eq-defX}), i.e.,
points satisfying the equality $\|A_{0}x\|=\|A_{1}x\|$. Figure~\ref{F:0b}
shows the trajectory with the maximum rate of increase of the Barabanov norm
$\|\cdot\|$. The construction of the next point of the trajectory $x_{n+1}$
depends on which of the sectors bounded by straight dash-dotted lines the
point $x_{n}$ belongs to; the type of matrix used for this, $A_{0}$ or
$A_{1}$, is indicated in the shaded areas of Fig.~\ref{F:0b}. The numerical
analysis performed shows that $\rho(\setA)\approx 1.098668$ and the index
sequence $\{\sigma_{n}\}$ of the extremal trajectory shown in
Fig.~\ref{F:0b} is of the form\footnote{In the theory of symbolic sequences,
it is customary to write the elements of the corresponding sequences in a row
without intermediate separators.}
\begin{equation}\label{E:indseq}
\{\sigma_{n}\}=10110110101101101101011011011010110110110101101101\ldots\,,
\end{equation}
where on sufficiently large segments (words) of the sequence $\{\sigma_{n}\}$
of length $10000$ the symbol $\boldsymbol{0}$ occurs with a frequency of
$\approx 0.364$ and the symbol $\boldsymbol{1}$  occurs with a frequency of
$\approx 0.636$. The index sequence $\{\sigma_{n}\}$ was constructed with the
program \texttt{\detokenize{barnorm_sturm.py}} by computing point iterations
using the angular function $\tilde{\varPhi}(\varphi)$, as shown in
Figs.~\ref{F:0c} and~\ref{F:01}. In these figures, the thick solid lines
denote sections of the graphs of the functions $\varPhi_{0}(\varphi)$ and
$\varPhi_{1}(\varphi )$ through which, according
to~\eqref{Eq-lxdef}--\eqref{E:angfun} the function $\tilde{\varPhi}(\varphi)$
is determined, and thin dashed lines mark those sections of the graphs of the
functions $\varPhi_{0}(\varphi)$ and $\varPhi_{1}(\varphi)$ that were
discarded in the definition of the function $\tilde{\varPhi}(\varphi)$.

In this case, the fact that both matrices $A_{0}$ and $A_{1}$ are nonnegative
and therefore leave the first quadrant in $\mathbb{R}^{2}$ invariant proved
to be of fundamental importance for the theoretical study of the structure of
extremal trajectories carried out in~\cite{Koz:INFOPROC05:e,
Koz:INFOPROC06:e}. The latter fact is reflected in the observation that the
angular function $\tilde{\varPhi}(\varphi)$ in this case maps the segment
$\left[0,\frac{\pi}{2}\right)$ into itself, which simplifies the study of its
trajectories $\{\tilde{\varphi}_{n}\}$. In particular, it was proved
in~\cite{Koz:INFOPROC05:e, Koz:INFOPROC06:e} that the restrictions of the
sets $\Omega_{0}$ and $\Omega_{1}$ in~\eqref{Eq-lxdef} on
$\left[0,\frac{\pi}{2}\right)$ are segments with a single common point
$\omega_{0}$. As numerical calculations show, in the case~\eqref{mycase} we
are considering, these sets are as follows:
\[\textstyle
\Omega_{0}\cap\left[0,\frac{\pi}{2}\right)=\left[\omega_{0},\frac{\pi}{2}\right),\quad
\Omega_{1}\cap\left[0,\frac{\pi}{2}\right)=\left[0,\omega_{0}\right).
\]
Accordingly, the function $\tilde{\varPhi}(\varphi)$ in our case has only one
discontinuity point $\omega_{0}$ on the segment
$\left[0,\frac{\pi}{2}\right)$, see Fig.~\ref{F:01}.

If we treat $\varphi\in\left[0,\frac{\pi}{2}\right)$ as an angular coordinate
on a circle of length $\frac{\pi}{2}$, we can take $\tilde{\varPhi}(\varphi)$
as an orientation-preserving mapping\footnote{The mapping $f(\varphi)$ of a
circle into itself is called \emph{orientation-preserving} if for any triple
of points $\varphi_{0},\varphi_{1},\varphi_{2}$ on the circle the order of
these points in going round the circle in any direction agrees with the order
of their images $f(\varphi_{0}),f(\varphi_{1}),f(\varphi_{2})$ when going
around the circle in the same direction.} of the corresponding circle into
itself, see~\cite{Koz:INFOPROC05:e, Koz:INFOPROC06:e} for details. In this
case, as shown in~\cite{Koz:INFOPROC05:e, Koz:INFOPROC06:e}, the index
sequence $\{\sigma_{n}\}$ of each trajectory $\tilde{\varphi}_{n}$ of the
mapping $\tilde{\varPhi}$ turns out to be either periodic or the so-called
\emph{Sturmian} (see, e.g.,~\cite[Ch.~6]{Fogg02}, \cite[Ch.~2]{Lothaire02}).

\begin{figure}[htbp!]
\centering \mbox{}\hfill\subcaptionbox{Rotation by the angle $\theta$ equal
to the length of the arc~$I_{1}$\label{F:rot-a}}
{\includegraphics*[height=0.31\textwidth]{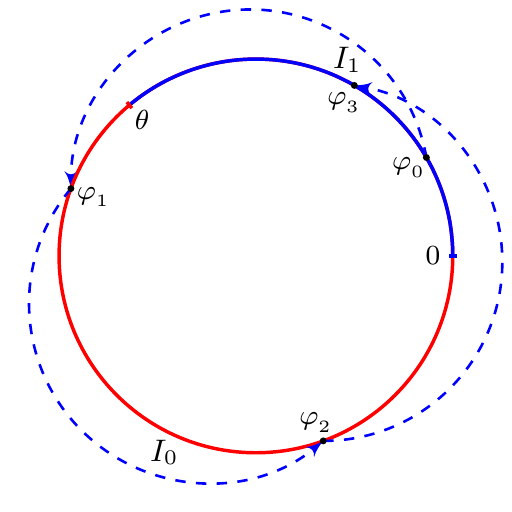}}
\hfill\subcaptionbox{Double rotation: each of the arcs $I_{0}$ and $I_{1}$
is rotated by its own angle\label{F:rot-b}}
{\includegraphics*[height=0.31\textwidth]{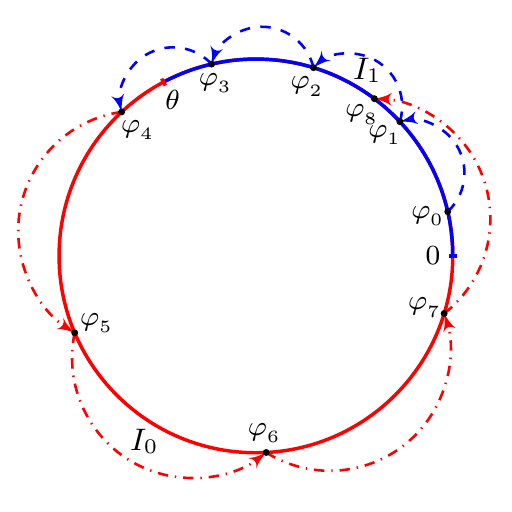}}
\hfill\mbox{}\caption{Regular and double rotation of a circle}\label{F:rot}
\end{figure}

Currently, there are several definitions of the Sturmian sequences. One of
the simplest analytic definitions states that a Sturmian sequence is an
integer sequence $\{\sigma_{n}\}$, which for all integer $n$ is defined by
the relations
\[
\sigma_{n}=[(n+1)\theta + \eta]-[n\theta +\eta],
\]
where $\theta\in(0,1)$ is an irrational number, $\eta\in\mathbb{R}$, and
$[\cdot]$ denotes the integer part of the number. The following equivalent
definition will be more useful for us: let $\{\varphi_{n}\}$ be a trajectory
running on a circle of length $1$ (or, equivalently, on the interval $[0,1)$)
through the rotation map
\begin{equation}\label{E:rotmap}
\varphi_{n+1} =\varphi_{n} +\theta\bmod1,
\end{equation}
where $\theta$ is an irrational number. We associate the index sequence
$\{\sigma_{n}\}$ with this trajectory and set (see Fig.~\ref{F:rot-a})
\begin{equation}\label{E:rotmap:i}
\sigma_{n}=\begin{cases}
0,&\text{if~}\varphi_{n}\in I_{0}:=[\theta,1),\\
1,&\text{if~}\varphi_{n}\in I_{1}:=[0,\theta).
\end{cases}
\end{equation}
Then the resulting sequence $\boldsymbol{\sigma}=\{\sigma_{n}\}$ (which is
not periodic due to the irrationality of $\theta$) is simply the so-called
Sturmian sequence formed by the symbol pair $\{0,1\}$ and the ``rotation
number'' $\theta$. In defining Sturmian sequences, we will later give up the
requirement that the number $\theta$ must be irrational. This will lead to
the fact that such generalized Sturmian sequences may turn out to be periodic.

We also note that the characteristic property of Sturmian sequences
$\boldsymbol{\sigma}$ with irrational $\theta$ is that they satisfy the
identity
\begin{equation}\label{E:subwordcompl}
	p(n,\boldsymbol{\sigma})\equiv n+1,
\end{equation}
where $p(n,\boldsymbol{\sigma})$ is the so-called
\emph{subword complexity function}, defined as the number of distinct words
of length $n$ in the sequence $\boldsymbol{\sigma}$, see,
e.g.,~\cite[Sec.~1.2.2]{Lothaire02}, \cite[Sec.~5.1.3]{Fogg02},
\cite[Sec.~6]{BK:BEATCS03}.

One of the most important properties of the Sturmian sequences is the following
fact~\cite[Lemma~6.1.3]{Fogg02}:

\begin{lemma}\label{L:sturmsym}
In any (generalized) Sturmian sequence consisting of two characters
$\{0,1\}$, exactly one of the symbol sequences (words) $\boldsymbol{00}$ or
$\boldsymbol{11}$ does not occur.
\end{lemma}

Lemma~\ref{L:sturmsym} becomes clear if we note that the points
$\{\varphi_{n}\}$ satisfying~\eqref{E:rotmap} cannot fall twice in succession
in the interval $[0,\theta )$ or $[\theta,1)$ which has the smallest length.

For illustration, note that \textbf{\emph{the Sturmian index
sequence~\eqref{E:indseq} does not contain the symbol sequence (word)
$\boldsymbol{00}$}.}

\section{A Pair of Matrices Similar to Plane Rotations}\label{S:rotmaps}
As mentioned in Section~\ref{S:sturmian-case}, previous
studies~\cite{LagWang:LAA95, BM:JAMS02, BTV:SIAMJMA03, Koz:INFOPROC05:e,
Koz:INFOPROC06:e} for matrix sets~\eqref{Eq-M1} consisting of nonnegative
matrices of a special form have provided some clarity on the structure of
index sequences that yield the maximal growth rate of matrix product norms.
For this reason, in this section we focus on considering less studied matrix
sets, namely sets consisting of matrices that may have negative elements. Our
goal is to give an example of matrix sets $\setA=\{A_{0},A_{1}\}$ consisting
of matrices of dimension ${2\times2}$, in which the sequences of indices
$\{\sigma_{i}\}$ maximizing $\|A_{\sigma_{n}}\cdots
A_{\sigma_{1}}A_{\sigma_{0}}x\|$, where $\|\cdot\|$ is a Barabanov norm,
\textbf{\emph{are not Sturmian}}! One of the simplest types of this kind of
matrix sets is the set of matrices $\setA=\{A_{0},A_{1}\}$, where
\begin{equation}\label{Eq-M2}
A_{0}=\left[\begin{array}{cr}
  \cos\theta_{0} & -\sin\theta_{0}\\
  \sin\theta_{0} & \cos\theta_{0}
\end{array}\right],\quad
A_{1}=\left[\begin{array}{rr}
  \cos\theta_{1} & -\lambda\sin\theta_{1}\\
  \frac{1}{\lambda}\sin\theta_{1} & \cos\theta_{1}
\end{array}\right].
\end{equation}

Further considerations in this section are based on computational experiments
and have no theoretical basis at this time. In this respect, this section
should be considered as a kind of set of questions (with accompanying
comments) for further research.

Consider the sets of matrices $\setA=\{A_{0},A_{1}\}$ defined by the following
parameters:
\begin{alignat*}{4}
&\text{Case 1}:\quad&\theta_{0}&=0.4,&\theta_{1}&=0.8,& \lambda&=0.75,\\
&\text{Case 2}:&\theta_{0}&\approx 0.6151,\quad&\theta_{1}&=0.8,\quad
&\lambda&=0.75,\\
&\text{Case 3}:&\theta_{0}&=0.7,&\theta_{1}&=0.8, &\lambda&=0.75.
\end{alignat*}

\begin{figure}[htbp!]
\centering \mbox{}\hfill\subcaptionbox{A Barabanov norm\label{F:1a}}
{\includegraphics*[height=0.31\textwidth]{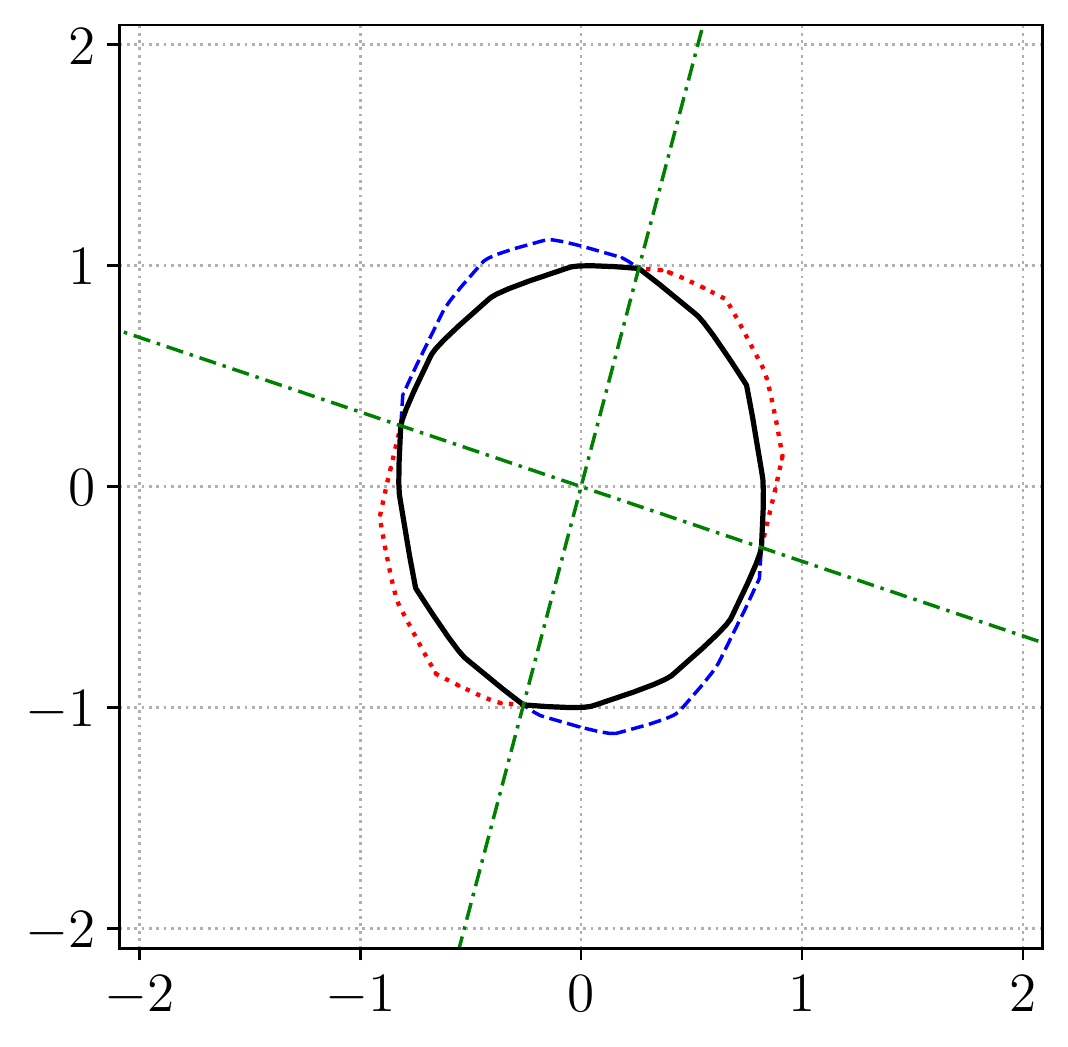}}
\hfill\subcaptionbox{A maximum growth rate trajectory\label{F:1b}}
{\includegraphics*[height=0.31\textwidth]{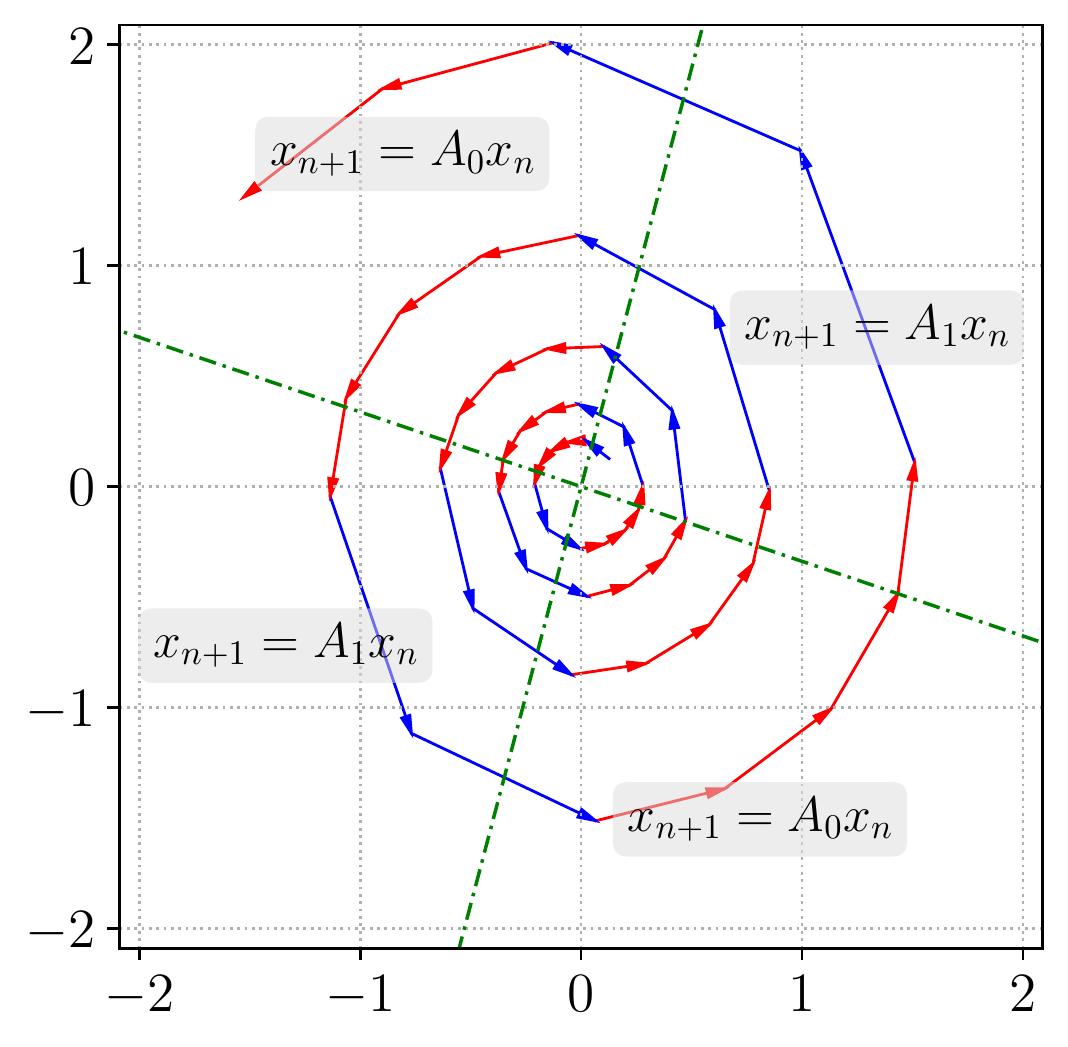}}
\hfill\subcaptionbox{The function $\tilde{\varPhi}(\varphi)$ does not
preserve orientation and has two discontinuity points\label{F:1c}}
{\includegraphics*[height=0.312\textwidth]{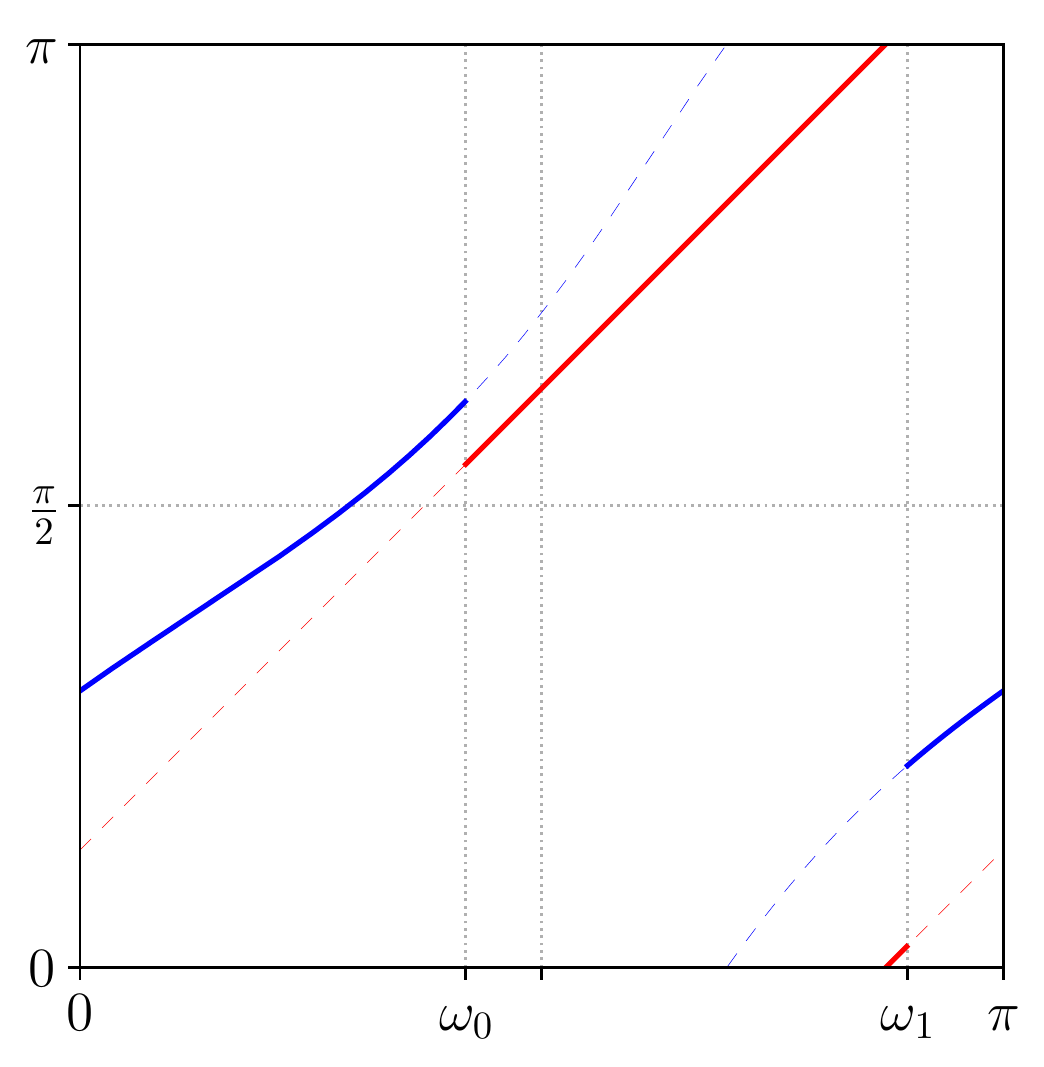}}
\hfill\mbox{} \caption{A Barabanov norm and the angular function:
Case~1}\label{F:1}
\end{figure}

\begin{figure}[htbp!]
\centering \mbox{}\hfill\subcaptionbox{A Barabanov norm\label{F:2a}}
{\includegraphics*[height=0.31\textwidth]{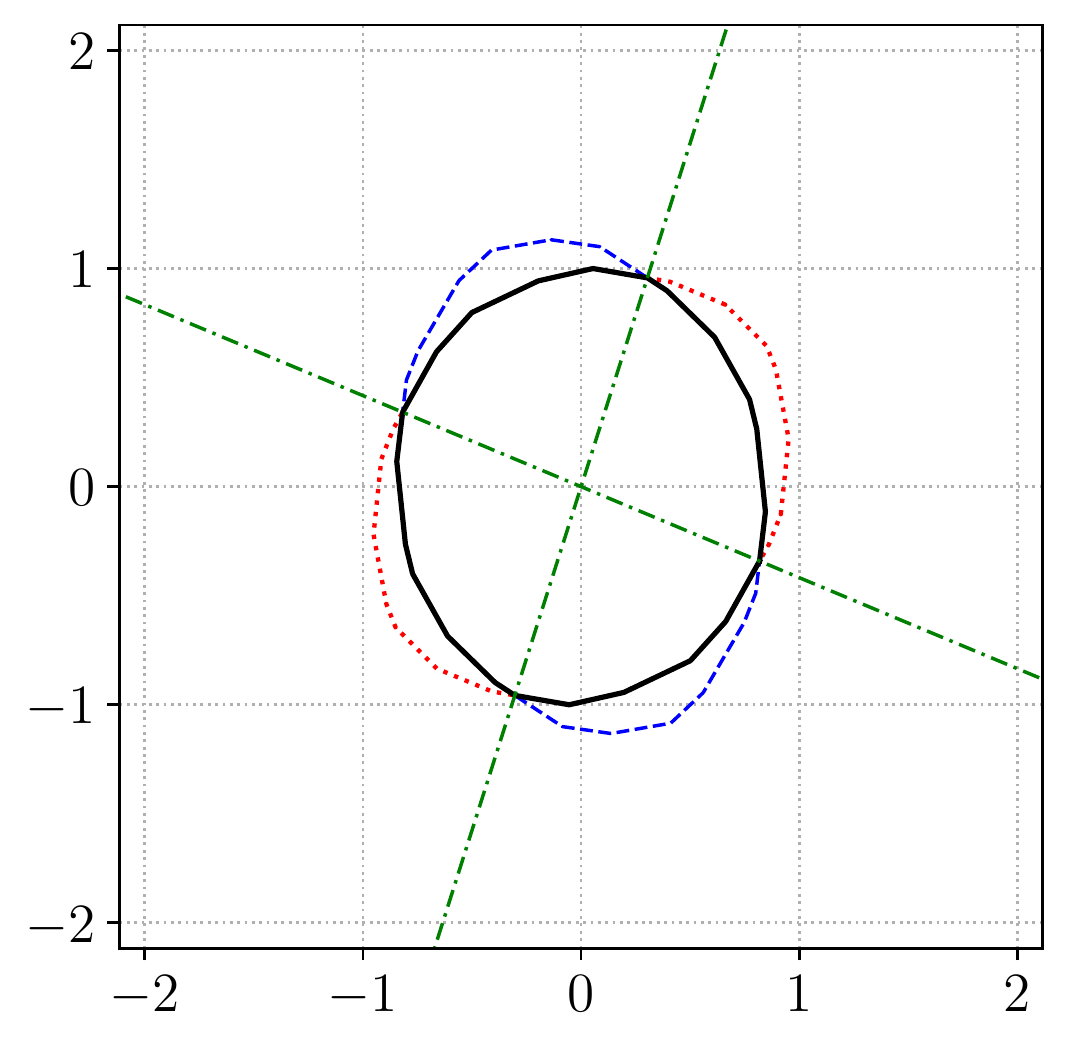}}
\hfill\subcaptionbox{A maximum growth rate trajectory\label{F:2b}}
{\includegraphics*[height=0.31\textwidth]{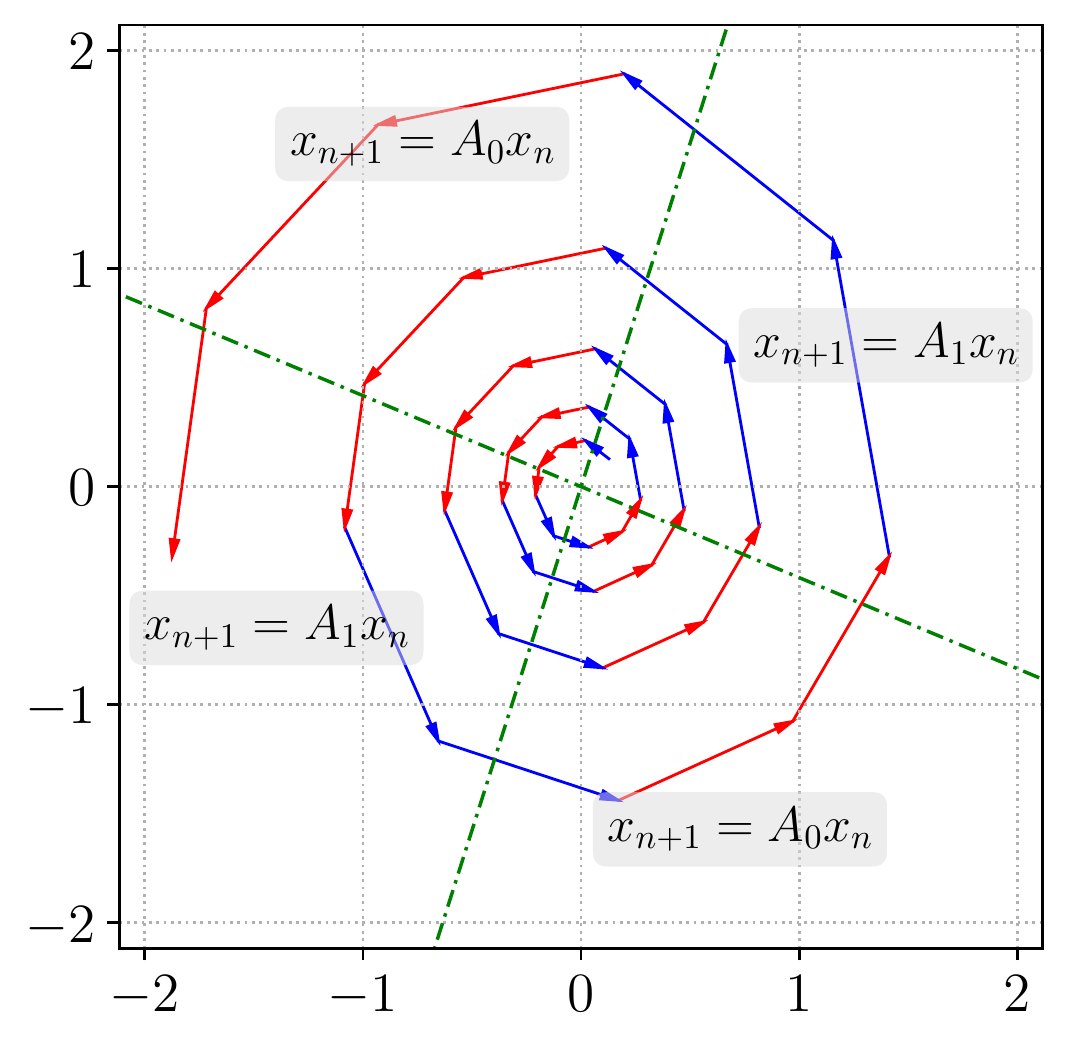}}
\hfill\subcaptionbox{The function $\tilde{\varPhi}(\varphi)$ preserves
orientation and has one discontinuity point\label{F:2c}}
{\includegraphics*[height=0.312\textwidth]{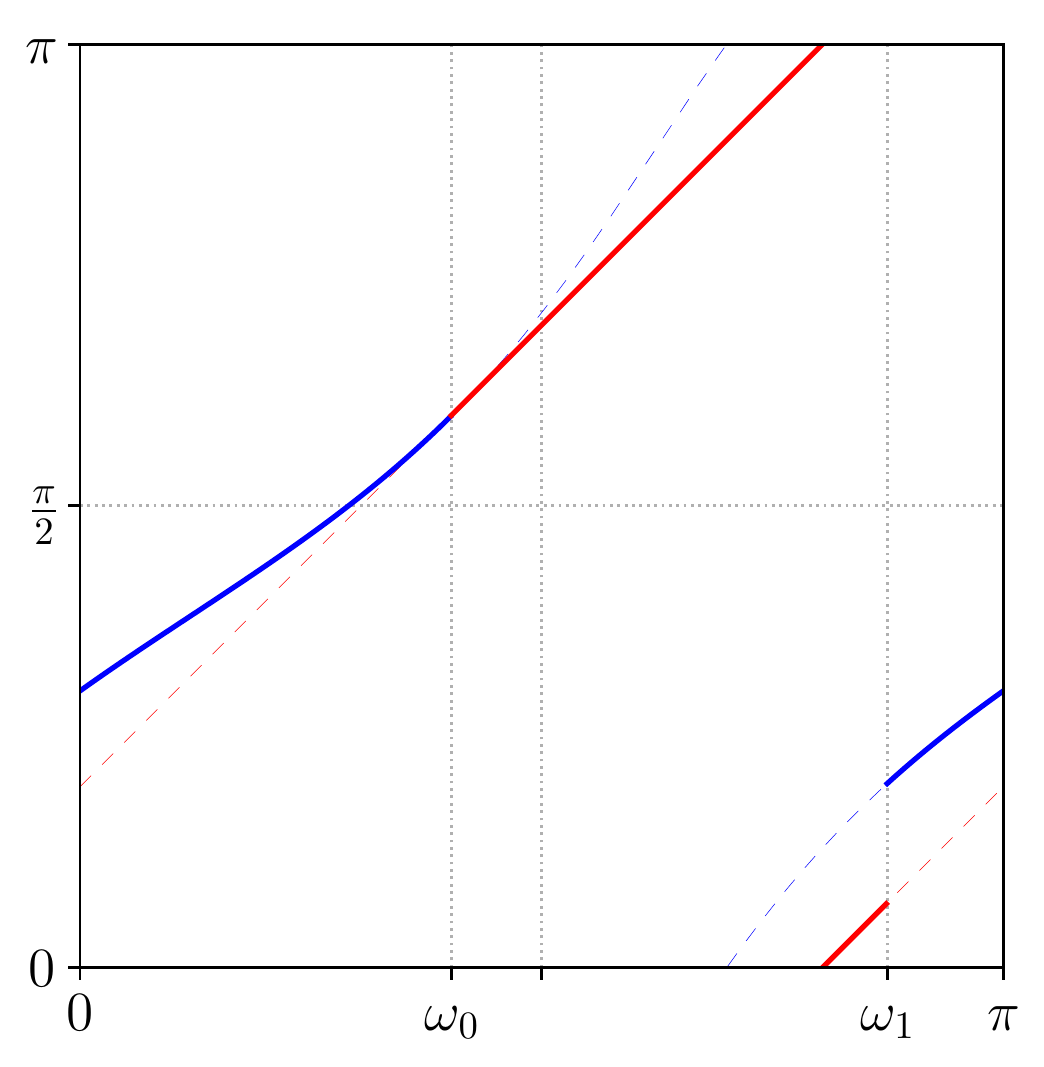}}
\hfill\mbox{} \caption{A Barabanov norm and the angular function:
Case~2}\label{F:2}
\end{figure}

\begin{figure}[htbp!]
\centering \mbox{}\hfill\subcaptionbox{A Barabanov norm\label{F:3a}}
{\includegraphics*[height=0.31\textwidth]{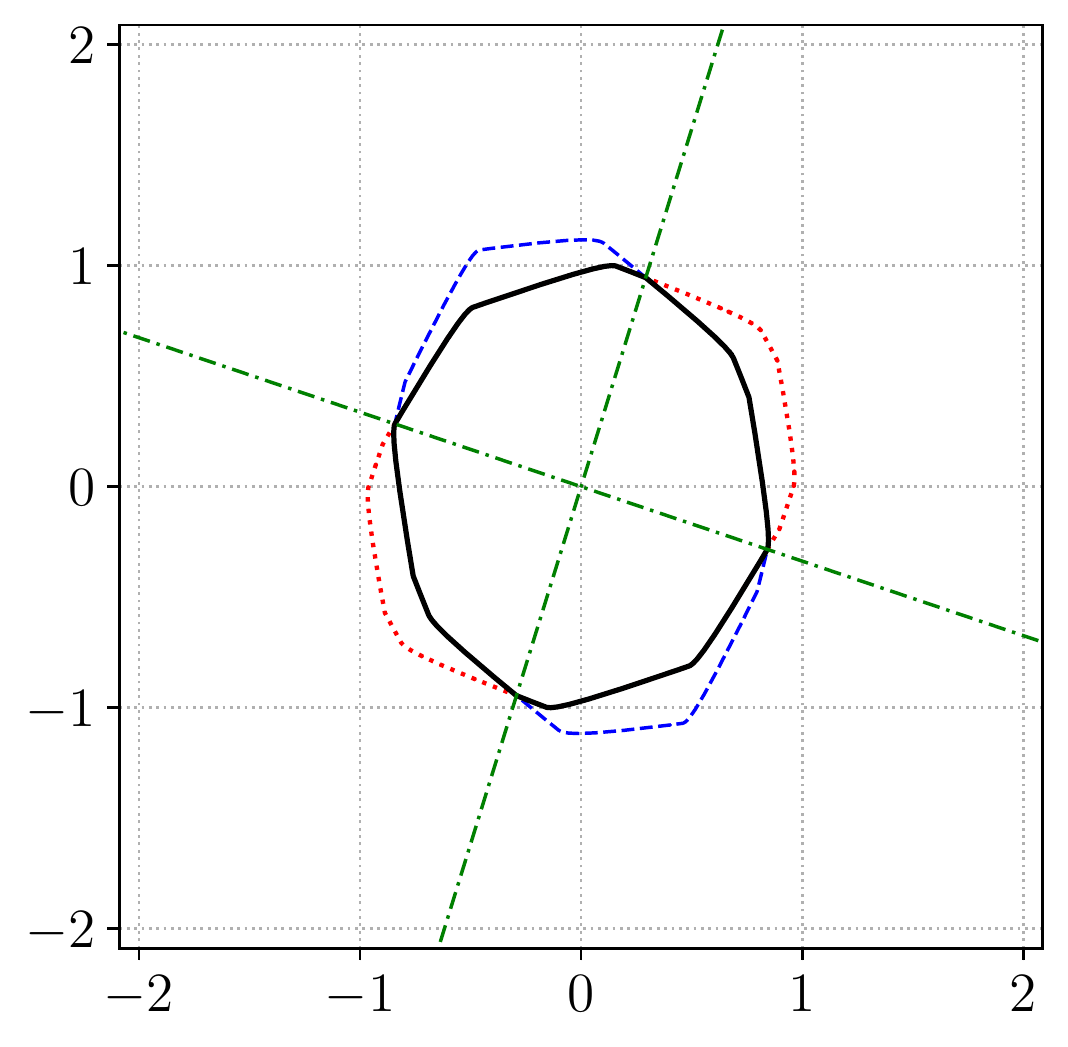}}
\hfill\subcaptionbox{A maximum growth rate trajectory\label{F:3b}}
{\includegraphics*[height=0.31\textwidth]{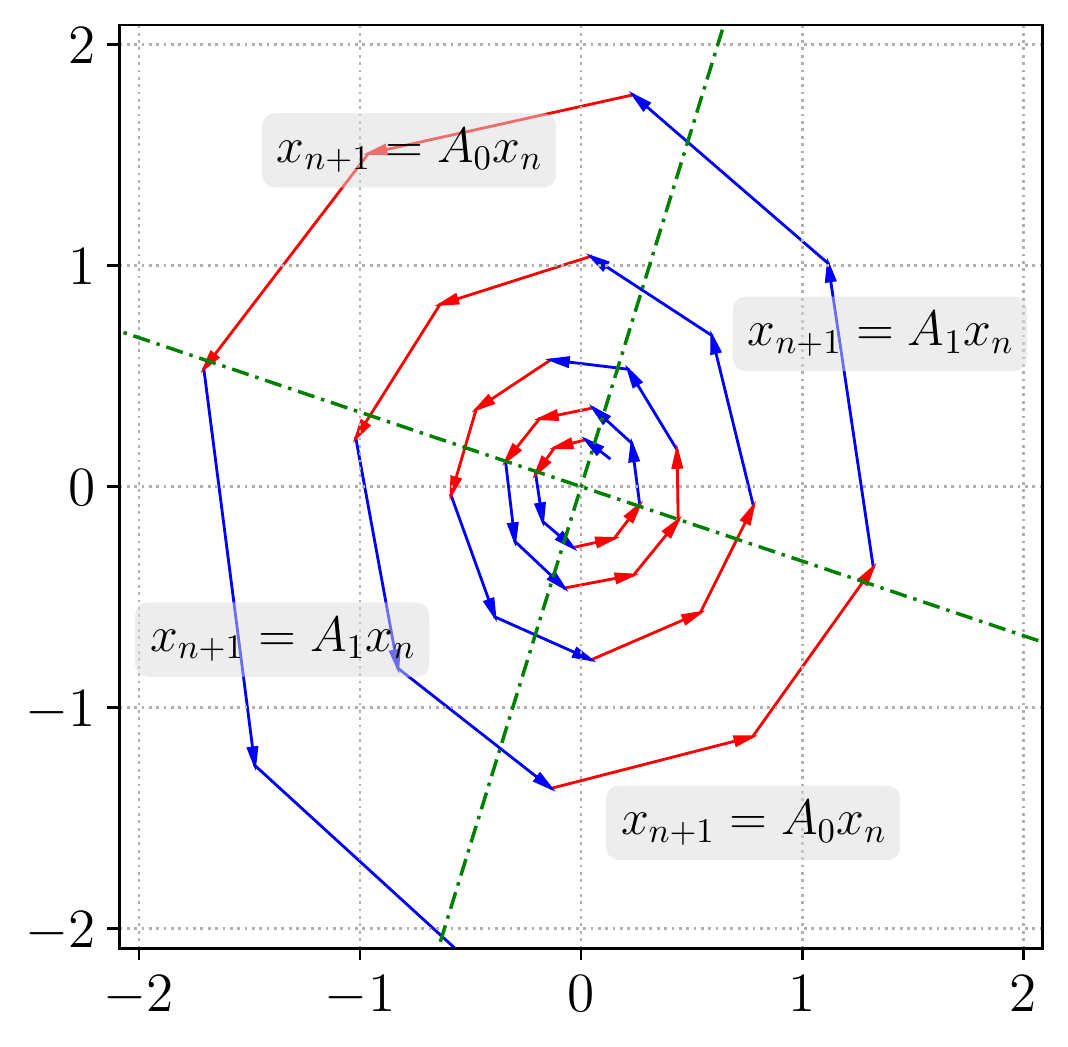}}
\hfill\subcaptionbox{The function $\tilde{\varPhi}(\varphi)$ preserves
orientation and has two discontinuity points\label{F:3c}}
{\includegraphics*[height=0.312\textwidth]{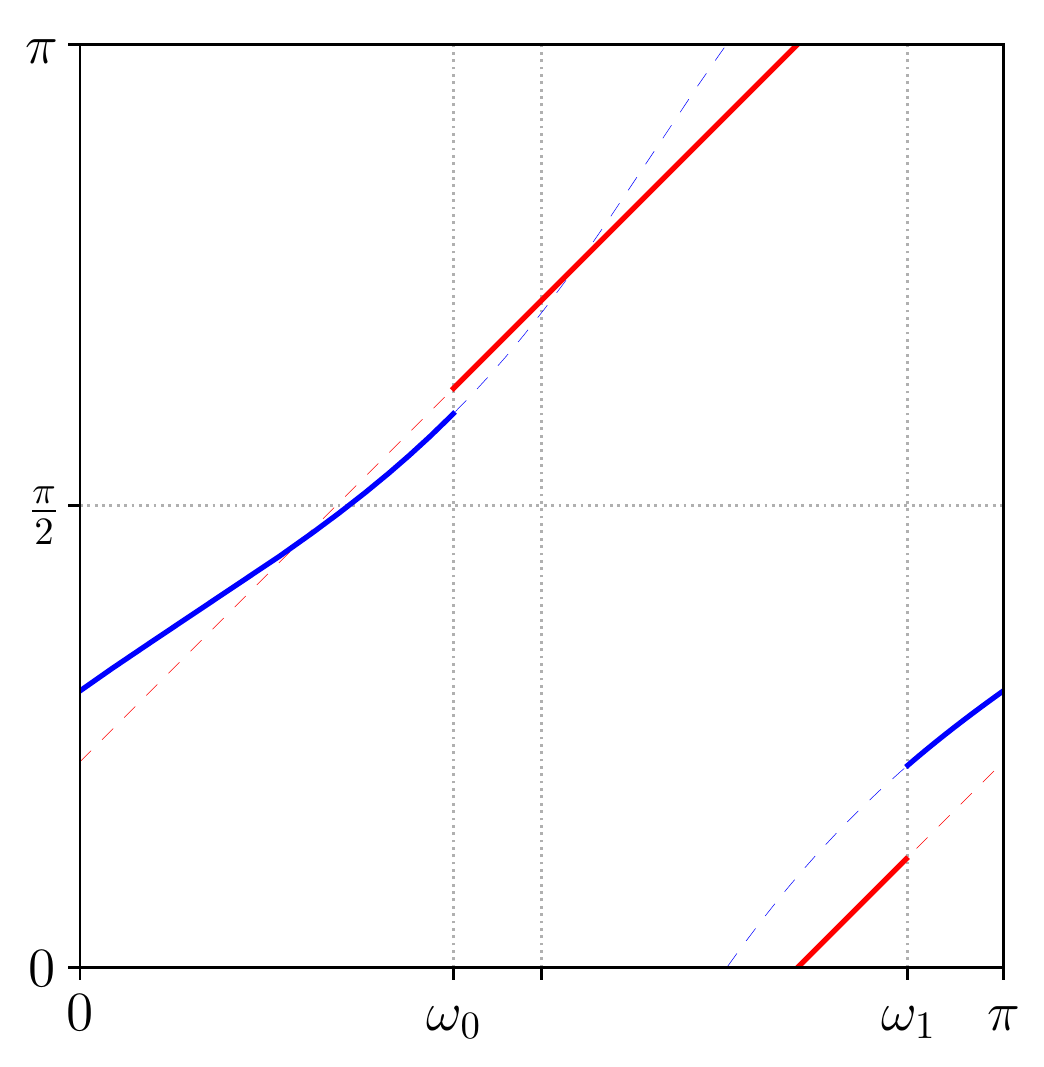}}
\hfill\mbox{} \caption{A Barabanov norm and the angular function:
Case~3}\label{F:3}
\end{figure}

In these cases, the software tools described in Section~\ref{S:numeric}
allowed not only to visualize approximately the shape of the unit sphere of
the Barabanov norm, but also to show examples of iterations
$x_{n+1}=A_{\sigma_{n}}x_{n}$ where the maximum growth rate of the Barabanov
norm of $\|x_{n}\|$ is reached. It also approximately finds the angular
function $\tilde{\varPhi}(\varphi)$ (see~\eqref{E:angfun}) of the matrix set
$\setA$, whose iterations allow the computation of the angular coordinates
$\tilde{\varphi}_{n}$ of the corresponding vectors $x_{n}$, without computing
their norms! The results of the corresponding numerical simulations are shown
in Figs.~\ref{F:1}, \ref{F:2} and \ref{F:3}. Since the meaning of the
notation in these figures repeats verbatim the explanations made for
Figs.~\ref{F:0a}, \ref{F:0b} and \ref{F:0c}, we do not present them here.

As can be seen in Figs.~\ref{F:1a}, \ref{F:2a}, and \ref{F:3a}, in all three
cases the set ${X_{0}\cap X_{1}}$ (see definitions in~\eqref{Eq-defX},
\eqref{Eq-maxmap}) turns out to be the union of two straight lines passing
through the origin (dash-dotted lines). In this context, the following question
arises.

\begin{question}\label{Q:1}\rm
For the case of a pair of nonnegative matrices~\eqref{Eq-M1}, the statement
that the part of the set $X_{0}\cap X_{1}$ passing through the first and
third quadrants is a straight line is strictly justified
in~\cite{Koz:INFOPROC05:e, Koz:INFOPROC06:e}. We are not aware of such a
proof for the sets of matrices~\eqref{Eq-M2} considered in this section. The
question arises: \textbf{Is the set $\boldsymbol{X_{0}\cap X_{1}}$ in this
case always the union of two straight lines? And why only of two?}\qed
\end{question}

From Figs.~\ref{F:1b}, \ref{F:2b}, and \ref{F:3b} it is evident (and it was
calculated with the program \texttt{\detokenize{barnorm_rot.py}} described in
Section~\ref{S:numeric}) that the index sequences of the trajectories with the
maximum growth rate in the Barabanov norm, for Cases~1--3, are as follows:
\begin{align*}
\{\sigma_{n}\}&=10000110000110000110001100001100001100001100001100\ldots\,,\\
\{\sigma_{n}\}&=10001100110001100110001100110001100110001100110001\ldots\,,\\
\{\sigma_{n}\}&=10011001100110001100110011001100110011001100011001\ldots\,.
\end{align*}
Calculations with the program \texttt{\detokenize{barnorm_rot.py}} have shown
that in sequences $\{\sigma_{n}\}$ of length 10000 the symbols
$\boldsymbol{0}$, $\boldsymbol{1}$, $\boldsymbol{00}$, $\boldsymbol{01}$,
$\boldsymbol{10}$ and $\boldsymbol{11}$ occur with the following frequencies:
\begin{center}
\begin{tabular}{lcccccc}
Symbols   & $\boldsymbol{0}$& $\boldsymbol{1}$& $\boldsymbol{00}$&
$\boldsymbol{01}$& $\boldsymbol{10}$& $\boldsymbol{11}$\\ \hline
Frequencies (Case 1)&0.655&0.345&0.483&0.172&0.172&0.172\\
Frequencies (Case 2)&0.555&0.445&0.333&0.222&0.222&0.222\\
Frequencies (Case 3)&0.517&0.483&0.276&0.241&0.241&0.241
\end{tabular}
\end{center}
From this, according to Lemma~\ref{L:sturmsym}, the following fundamentally
important conclusion follows.

\begin{mainclaim}
\makeatletter\def\@currentlabelname{Main Claim}\makeatother%
\label{P:main}%
In the cases of pairs of matrices~\eqref{Eq-M2}, similar to plane rotations,
the index sequences of trajectories with the maximal growth rate in the
Barabanov norm \textbf{are not (generalized) Sturmian}.
\end{mainclaim}

Let us turn to a more detailed analysis of the obtained results of numerical
simulation.

\paragraph{Case 1.} As can be seen from Fig.~\ref{F:1c}, in Case~1 the
angular function $\tilde{\varPhi}(\varphi)$, considered as a mapping of the
circle into itself, \textbf{does not preserve orientation}. In particular,
this situation resembles the behavior of the so-called \emph{double
rotations}~\cite{SIA:DCDS05, Clack:13, AFHS:ArXiv21, Kryzhevich:MMNP21},
defined by the equation
\begin{equation}\label{E:doublerot}
f_{\theta_{1},\theta_{2},\theta}(\varphi) =
\begin{cases}
\varphi + \theta_{1}\bmod1& \text{for~}x\in [0, \theta),\\
\varphi + \theta_{2}\bmod1& \text{for~}x\in [\theta, 1),
\end{cases}
\end{equation}
where $(\theta_{1}, \theta_{2}, \theta)\in [0, 1)\times [0, 1)\times [0, 1)$
are some parameters (see the example in Fig.~\ref{F:rot-b}).

The mappings $\tilde{\varPhi}(\varphi)$ and
$f_{\theta_{1},\theta_{2},\theta}(\varphi)$ are related by the fact that both
are not continuous and do not preserve orientation on a circle but their
action is determined by some ``rotations'' on two continuity intervals of the
mappings. The difference between the angular function
$\tilde{\varPhi}(\varphi)$ from the double rotations of the circle
$f_{\theta_{1},\theta_{2},\theta}(\varphi)$ is that for the first of these
mappings, the rotation angles are not constant, while in the case of the
mapping $f_{\theta_{1},\theta_{2},\theta}(\varphi)$, the rotation angles are
constant for each of the intervals $[0, \theta)$ and $[\theta, 1)$. Double
rotations of the circle are thus somewhat easier to study, and some progress
has been made recently in their analysis~\cite{SIA:DCDS05, Clack:13,
AFHS:ArXiv21}.

\begin{question}\label{Q:2}\rm
Is it possible (by analogy with the case of the angular function
$\tilde{\varPhi}(\varphi)$ for the set of nonnegative matrices~\eqref{Eq-M1})
for the angular function $\tilde{\varPhi}(\varphi)$, arising in Case~1,
select parameters $\theta_{1}$, $\theta_{2}$, $\theta$ such that the
corresponding index sequences for the mapping $\tilde{\varPhi}(\varphi)$
would match the index sequences for the mapping
$f(\theta_{1},\theta_{2},\theta)(\varphi)$?\qed
\end{question}

A positive answer to this question does not seem very likely. As a first step
to clarify the situation, it would be possible to compare the frequencies of
occurrence of the symbols $\boldsymbol{0}$ and $\boldsymbol{1}$, as well as
groups of consecutive identical symbols $\boldsymbol{00\ldots0}$ and
$\boldsymbol{11\ldots1}$ in index sequences for the mappings
$\tilde{\varPhi}(\varphi)$ and $f(\theta_{1},\theta_{2},\theta)(\varphi)$.
Perhaps a negative answer to Question~\ref{Q:2} could have been obtained
already at this stage.

\begin{question}\label{Q:3}\rm
Again by analogy with the case of the angular function
$\tilde{\varPhi}(\varphi)$ for the set of matrices~\eqref{Eq-M1}: For the
angular function $\tilde{\varPhi}(\varphi)$ occurring in Case~1, are there
limit frequencies for the occurrence of the symbols $\boldsymbol{0}$ and
$\boldsymbol{1}$ in index sequences $\{\sigma_{n}\}$? If the answer is yes,
do these frequencies depend on a particular index sequence or not (as in the
case of the angular function $\tilde{\varPhi}(\varphi)$ for the matrix
set~\eqref{Eq-M1})?\qed
\end{question}

\begin{question}\label{Q:4}\rm
If Question~\ref{Q:3} is answered in the affirmative, then are the limiting
frequencies of occurrence of the symbols $\boldsymbol{0}$ and
$\boldsymbol{1}$ in the index sequences $\{\sigma_{n}\}$ for the angular
function $\tilde{\varPhi}(\varphi)$ dependent on a particular index sequence
or not (as in the case of the angular function $\tilde{\varPhi}(\varphi)$ for
a set of matrices~\eqref{Eq-M1})?\qed
\end{question}

In answering Question~\ref{Q:4}, it might be useful to refer to the theory of
circle mappings~\cite{Misiurewicz:ETDS86, ALM:N90, AM:AMUC96, AM:AMUC00,
Misiurewicz07}, both continuous and discontinuous, which are not
orientation-preserving. Unfortunately, the lack of the orientation preserving
property makes the analysis of the mappings of the circle much more
difficult. Instead of characterizing the ``mean rotation angle'' in such
mappings by the so-called ``rotation number'' (a standard tool in the theory
of orientation-preserving circle mappings), the concept of a ``rotation
interval'' emerges in non-orientation-preserving circle
mappings~\cite{Misiurewicz:ETDS86, ALM:N90, AM:AMUC00, Misiurewicz07}. The
latter fact may be crucial in answering the question whether the frequency of
occurrence of symbols in an index sequence depends on that sequence. Note,
however, that in our numerical experiments the dependence of the frequency
characteristics on the trajectory was not observed.

\paragraph{Cases~2 and~3.} In these cases the angular function
$\tilde{\varPhi}(\varphi)$, considered as a mapping of the circle into
itself, \textbf{preserves orientation}. In Case~2 it has one discontinuity
point and in Case~3 two, see Figs.~\ref{F:2c} and~\ref{F:3c}.
Therefore~\cite{Koz:CDC05:e, Koz:INFOPROC05:e, Koz:INFOPROC06:e}, for the
mapping $\tilde{\varPhi}(\varphi)$ in these cases, the so-called
\emph{rotation number} $\varkappa(\tilde{\varPhi})$ is defined, which
characterizes the ``mean rotation angle'' performed by this mapping.

In the case where $\tilde{\varPhi}(\varphi)$ is an angular
function~\eqref{Eq-maxmanlx} generated by the pair of nonnegative
matrices~\eqref{Eq-M1}, the rotation number $\varkappa(\tilde{\varPhi})$
coincides with the frequency of hitting the trajectory elements
\[
x_{n+1}=A_{\sigma_{n}}x_{n},\qquad n=0,1,\ldots~,
\]
into the set $X_{0}$, and hence with the frequency of occurrence of the
symbol $\boldsymbol{0}$ in the corresponding index sequence.

\begin{question}\label{Q:5}\rm
In Cases~2 and~3, as mentioned above, a rotation number is also defined for
the mapping $\tilde{\varPhi}(\varphi)$. However, whether this implies the
existence of limiting frequencies for the occurrence of the symbols
$\boldsymbol{0}$ and $\boldsymbol{1}$ in the index sequences $\{\sigma_{n}\}$
for the angular function $\tilde{\varPhi}(\varphi)$ remains unclear!\qed
\end{question}

The behavior of the mapping $\tilde{\varPhi}(\varphi)$, considered as a
mapping of a circle, resembles the behavior of the mapping of a
circle~\eqref{E:rotmap} in Cases~2 and~3:
\begin{equation}\label{E:rotmap2}
\varphi_{n+1} = \varphi_{n} +\theta\bmod1,
\end{equation}
with the difference that this time the index sequence $\{\sigma_{n}\}$ is not
calculated using formula~\eqref{E:rotmap:i}, but as follows:
\begin{equation}\label{E:rotmap2:i}
\sigma_{n}=\begin{cases}
0,&\text{if~}\varphi_{n}\in I_{0}:=[\theta_{0},1),\\
1,&\text{if~}\varphi_{n}\in I_{1}:=[0,\theta_{0}),
\end{cases}
\end{equation}
where the length of the interval $I_{1}$ is generally different from the
angle of rotation: $\theta_{0}\neq\theta$. As shown in~\cite{BV:TCS02}, the
behavior of the index sequence~\eqref{E:rotmap2:i} can be expressed in terms
of a pair of Sturmian sequences generated by rotation through angle $\theta$,
but in a rather complex way.

A question similar to that of Question~\ref{Q:2} may be asked here.

\begin{question}\label{Q:6}\rm
For the angular function $\tilde{\varPhi}(\varphi)$ occurring in Cases~2
and~3, is it possible to choose the parameters $\theta_{0}$ and $\theta$ such
that the corresponding index sequences for $\tilde{\varPhi}(\varphi)$
coincide with the index sequences~\eqref{E:rotmap2:i} for the rotation
mapping~\eqref{E:rotmap2}?\qed
\end{question}

As in the case of Question~\ref{Q:2}, a positive answer to this question does
not seem very likely. Here, as a first step to clarify the situation, one
might compare the frequencies of occurrence of the symbols $\boldsymbol{0}$
and $\boldsymbol{1}$ as well as groups of consecutive identical symbols
$\boldsymbol{00\ldots0}$ and $\boldsymbol{11\ldots1}$ in index sequences for
the mappings $\tilde{\varPhi}(\varphi)$
and~\eqref{E:rotmap2}--\eqref{E:rotmap2:i}. Note, however, that in this case
(unlike the situation described in the discussion of Question~\ref{Q:4}) the
computation of the frequency characteristics of the index
sequences~\eqref{E:rotmap2:i} can be performed theoretically, which is
possible may simplify the research.

\section{Methods and Tools for Numerical Modeling}\label{S:numeric}

In this work, an approximate construction of Barabanov norms of matrix sets
(and visualization of their unit spheres as well as the trajectories with the
maximum growth rate in the Barabanov norm) was performed using the programs
\texttt{\detokenize{barnorm_sturm.py}} and
\texttt{\detokenize{barnorm_rot.py}}, which are available for download from
the website \url{https://github.com/kozyakin/barnorm}. These programs use a
small modification of the max-relaxation algorithm for the iterative
construction of Barabanov norms, which can be found in \cite{Koz:DCDSB10,
Koz:JDEA11}. The modification compared to the software implementation of the
corresponding algorithms described in~\cite{Koz:ArXiv10-1} was that convex
centrally symmetric polygons were chosen as unit spheres of norms
approximating the Barabanov norm. The advantage of this approach over the
approach of~\cite{Koz:ArXiv10-1} is that when linear transformations are
applied, the unit spheres of the norms $\|A_{0}x\|$ and $\|A_{1}x\|$ are
again convex centrally symmetric polygons. Using the library \texttt{shapely}
of the language \texttt{Python}, this allows the iterative computation of the
norm $\max\{\|A_{0}x\|, \|A_{1}x\|\}$ without loss of precision for each
iteration.

The programs \texttt{\detokenize{barnorm_sturm.py}} and
\texttt{\detokenize{barnorm_rot.py}} differ from each other only in the
specification of the matrices $A_{0}$ and $A_{1}$ and in the amount of
graphical data displayed. Both programs are implemented in the Python
programming language of versions 3.8--3.10 of the \texttt{Miniconda3}
(\texttt{Anaconda3}) distribution.

For the convenience of the reader, a listing of the program
\texttt{\detokenize{barnorm_rot.py}} is provided in Appendix~\ref{A:PyCode}.

\section{Conclusion}\label{S:conclusion}

The paper presents the results of a numerical simulation of the fastest
growing (in the Barabanov norm) trajectories generated by sets of
${2\times2}$ matrices. The results obtained indicate that in certain
situations the maximum growth rate can be achieved on trajectories with
non-Sturmian sequences of indices, which makes these situations fundamentally
different from most theoretical studies carried out so far in the theory of
joint/generalized spectral radius.

Section~\ref{S:rotmaps} presents the results of the numerical simulations
performed in this paper and formulates a number of open questions. In
particular, it would be interesting to compare the complexity functions
$p(n,\boldsymbol{\sigma})$ of the index sequences
$\boldsymbol{\sigma}=\{\sigma_{n}\}$ of the mappings
$\tilde{\varPhi}(\varphi)$ with the complexity functions of the index
sequences of the ``test'' circle mappings~\eqref{E:doublerot}
and~\eqref{E:rotmap2}--\eqref{E:rotmap2:i}. This question seems all the more
interesting because, unlike Sturmian sequences for which the growth rate of
the function $p(n,\boldsymbol{\sigma})$ is linear
(see~\eqref{E:subwordcompl}), for sequences generated by the double rotation
mappings~\eqref{E:rotmap2}--\eqref{E:rotmap2:i}, the growth rate of the
complexity function can be superlinear: $p(n,\boldsymbol{\sigma})\sim
n^{\gamma}$, where $\gamma > 1$~\cite{Clack:13}. An example of such a growth
rate of the norms of matrix products is described in~\cite{HMS:MPCPS13}.

Since the bulk of the results presented above are numerical in nature, this
paper should not be considered a full-fledged theoretical study but rather a
plan for further research on the subject.

\section*{Acknowledgment}
The author thanks Aljo{\v{s}}a Peperko for pointing out some obscurities in
the presentation.


\clearpage
\appendix
\section{Program for Calculating the Barabanov Norm}\label{A:PyCode}
The code below is written in the \texttt{Python} programming language
versions 3.8--3.10 of the \texttt{Miniconda3} (\texttt{Anaconda3})
distribution. This and some other related scripts for calculating the
Barabanov norm can be downloaded from the website
\url{https://github.com/kozyakin/barnorm}.
\medskip

\lstinputlisting[caption={Python code \texttt{\detokenize{barnorm_rot.py}}
for computing the Barabanov norm of a pair of matrices and the angular
function of the iterations at which the joint spectral radius is reached},
label=L:code]{\detokenize{barnorm_rot.py}}

\end{document}